\documentclass[aop,seceqn,MSNbibl,citesort,dvips]{arximspdf}
\usepackage{mathbh}
\usepackage{graphicx}

\doi{10.1214/10-AOP587}
\volume{39}
\issue{4}
\pubyear{2011}
\firstpage{1544}
\lastpage{1590}

\makeatletter

\newtheorem{theorem}{Theorem}[section]
\newtheorem{prop}[theorem]{Proposition}
\newtheorem{lem}[theorem]{Lemma}
\newproclaim{defi}[theorem]{Definition}
\newproclaim{rem}[theorem]{Remark}

\newcommand{\eset}{\varnothing}
\newcommand{\si}{\sigma}
\newcommand{\wt}{\widetilde}
\newcommand{\ind}{\mathbh{1}}
\newcommand{\sign}{\operatorname{sign}}
\newcommand{\tr}{\operatorname{tr}}
\newcommand{\pwit}{\operatorname{PWIT}}
\newcommand{\bv}{\mathbf{v}}
\newcommand{\bw}{\mathbf{w}}
\newcommand{\tbT}{\wt{\mathbf{T}}}

\newcommand{\cA}{\mathcal{A}}
\newcommand{\cD}{\mathcal{D}}
\newcommand{\cE}{\mathcal{E}}
\newcommand{\cK}{\mathcal{K}}
\newcommand{\cL}{\mathcal{L}}
\newcommand{\cM}{\mathcal{M}}
\newcommand{\cS}{\mathcal{S}}
\newcommand{\cT}{\mathcal{T}}
\newcommand{\cV}{\mathcal{V}}
\newcommand{\cY}{\mathcal{Y}}
\newcommand{\cZ}{\mathcal{Z}}
\newcommand{\bA}{\mathbf{A}}
\newcommand{\bD}{\mathbf{D}}
\newcommand{\bK}{\mathbf{K}}
\newcommand{\bN}{\mathbf{N}}
\newcommand{\bS}{\mathbf{S}}
\newcommand{\bT}{\mathbf{T}}
\newcommand{\bU}{\mathbf{U}}
\newcommand{\dP}{\mathbb{P}}

\makeatother

\begin{document}
\begin{frontmatter}

\title{Spectrum of large random reversible Markov~chains:
Heavy-tailed weights on~the~complete graph}
\runtitle{Heavy-tailed weights on the complete graph}

\begin{aug}
\author[A]{\fnms{Charles} \snm{Bordenave}\ead[label=e1]{charles.bordenave(at)math.univ-toulouse.fr}\ead[label=u1,url]{http://www.math.univ-toulouse.fr/\textasciitilde bordenave/}},
\author[C]{\fnms{Pietro} \snm{Caputo}\thanksref{t1}\ead[label=e2]{caputo(at)mat.uniroma3.it}\ead[label=u2,url]{http://www.mat.uniroma3.it/users/caputo/}} and
\author[B]{\fnms{Djalil} \snm{Chafa\"{i}}\corref{}\ead[label=e3]{djalil(at)chafai.net}\ead[label=u3,url]{http://djalil.chafai.net/}}

\runauthor{C. Bordenave, P. Caputo and D. Chafa\"{i}}
\affiliation{Universit\'{e} de Toulouse III, Universit\`{a} Roma
Tre and~Universit\'{e}~Paris-Est~Marne-la-Vall\'{e}e}

\address[A]{C. Bordenave\\
IMT UMR5219 CNRS\\
Universit\'{e} de Toulouse III\\
France\\
\printead{e1}\\
\printead{u1}}
\address[B]{D. Chafa{\"{i}}\\
LAMA UMR8050 CNRS\\
Universit\'{e} Paris-Est Marne-la-Vall\'{e}e\\
France\\
\printead{e3}\\
\printead{u3}}
\address[C]{P. Caputo\\
Dipartimento di Matematica\\
Universit\`{a} Roma Tre\\
Italy\\
\printead{e2}\\
\printead{u2}}
\end{aug}

\thankstext{t1}{Supported in part by NSF Grant
DMS-03-01795 and by the European Research Council through the ``Advanced
Grant'' PTRELSS 228032.}

\received{\smonth{4} \syear{2009}}
\revised{\smonth{5} \syear{2010}}

%
\begin{abstract}
We consider the random reversible Markov kernel $K$ obtained by assigning
i.i.d. nonnegative weights to the edges of the complete graph over $n$
vertices and normalizing by the corresponding row sum. The weights are
assumed to be in the domain of attraction of an $\alpha$-stable law,
$\alpha\in(0,2)$. When $1\leq\alpha<2$, we show that for a suitable
regularly varying sequence $\kappa_n$ of index $1-1/\alpha$, the limiting
spectral distribution $\mu_\alpha$ of $\kappa_nK$ coincides with the
one of
the random symmetric matrix of the un-normalized weights (L\'{e}vy
matrix with
i.i.d. entries). In contrast, when $0<\alpha<1$, we show that the empirical
spectral distribution of $K$ converges without rescaling to a nontrivial
law $\wt\mu_\alpha$ supported on $[-1,1]$, whose moments are the return
probabilities of the random walk on the Poisson weighted infinite tree
(PWIT) introduced by Aldous. The limiting spectral distributions are given
by the expected value of the random spectral measure at the root of suitable
self-adjoint operators defined on the PWIT. This characterization is used
together with recursive relations on the tree to derive some
properties of
$\mu_\alpha$ and $\wt\mu_\alpha$. We also study the limiting behavior
of the
invariant probability measure of $K$.
\end{abstract}

%

\end{frontmatter}

\section{Introduction}
Let $G_n=(V_n,E_n)$ denote the complete graph with vertex set
$V_n=\{1,\ldots,n\}$, and edge set $E_n=\{ \{i,j\} , 1\leq i,j\leq
n\}$,
including loops $\{i,i\}$, $1\leq i\leq n$. Assign a nonnegative random
weight (or conductance) $U_{i,j}=U_{j,i}$ to each edge $\{i,j\}\in
E_n$, and
assume that the symmetric weights $\bU=\{U_{i,j}; \{i,j\}\in E_n\}$ are i.i.d. with
common law $\mathcal{L}$ independent of $n$. This defines a random
network, or
weighted graph, denoted $(G_n,\bU)$. Next, consider the random walk on
$(G_n,\bU)$ defined by the transition probabilities
%
%
\begin{equation}\label{srw}
K_{i,j} := \frac{U_{i,j}}{\rho_i} \qquad
\mbox{with }
\rho_i:=\sum_{j=1}^n U_{i,j}.
\end{equation}
The Markov kernel $K$ is \textit{reversible} with respect to the measure
$\rho=\sum_{i\in V_n}\rho_i\delta_i$ in that
\[
\rho_i K_{i,j} = \rho_j K_{j,i}
\]
for all $i,j\in V_n$. Note that we have not assumed that $\cL$ has no
atom at
$0$. If $\rho_i=0$ for some $i$, then for that index $i$ we set
$K_{i,j}=\delta_{i,j}$, $1\leq j\leq n$. However, as soon as $\cL$ is not
concentrated at $0$ then almost surely, for all $n$ sufficiently large,
$\rho_i>0$ for all $1\leq i\leq n$, $K$ is irreducible and aperiodic and
$\rho$ is its unique invariant measure, up to normalization (see,
e.g.,~\cite{bordenave-caputo-chafai}).

For any square $n\times n$ matrix $M$ with eigenvalues
$\lambda_1(M),\ldots,\lambda_n(M)$, the \textit{Empirical Spectral
Distribution}
(ESD) is the discrete probability measure with at most $n$ atoms
defined by
\[
\mu_M:=\frac1n\sum_{j=1}^n \delta_{\lambda_j(M)}.
\]
All matrices $M$ to be considered in this work have real spectrum, and the
eigenvalues will be labeled in such a way that
$\lambda_n(M)\leq\cdots\leq\lambda_1(M)$.

Note that $K$ defines a square $n\times n$ random Markov matrix
whose entries are not independent due the normalizing sums $\rho_i$. By
reversibility, $K$ is self-adjoint in $L^2(\rho)$ and its spectrum
$\si(K)$
is real. Moreover, $\si(K)\subset[-1,+1]$, and $1\in\si(K)$.
Since $K$ is Markov, its ESD $\mu_K$ carries further probabilistic content.
Namely, for any $\ell\in\mathbb{N}$, if $p_\ell(i)$ denotes the probability
that the random walk on $(G_n,\bU)$ started at $i$ returns to $i$
after $\ell$
steps, then the $\ell$th moment of $\mu_K$ satisfies
%
%
\begin{equation}\label{moms}
\int_{-1}^{+1} x^\ell\mu_K(dx) %
= \frac1n\tr(K^\ell) %
= \frac1n\sum_{i\in V} p_\ell(i).
\end{equation}
%

\subsubsection*{Convergence of the ESD}
The asymptotic behavior of $\mu_K$ as $n\to\infty$ depends 
strongly on the tail of $\cL$ at infinity. When $\cL$ has finite mean
$\int_0^\infty x \cL(dx) = m$ we set $m=1$. This is no loss of generality
since $K$ is invariant under the dilation $t\to t U_{i,j}$. If $\cL$
has a
finite second moment we write $\si^2 = \int_0^\infty(x-1)^2 \cL
(dx)$ for the
variance.

The following result, from~\cite{bordenave-caputo-chafai}, states that if
$0<\si^2<\infty$, then the bulk of the spectrum of $\sqrt{n}K$
behaves, when
$n\to\infty$, as if we had truly i.i.d. entries (Wigner matrix).
Without loss
of generality, we assume that the weights $\bU$ come from the
truncation of a
unique infinite table $(U_{i,j})_{{i,j}\geq1}$ of i.i.d. random
variables of law
$\cL$. This gives a meaning to the almost sure (a.s.) convergence of
$\mu_{\sqrt n K}$. The symbol $\stackrel{w}{\to}$ denotes weak
convergence of measures with respect to continuous bounded functions. Note
that $\lambda_1(\sqrt{n}K)=\sqrt{n}\to\infty$.
\begin{theorem}[(Wigner-like behavior)]\label{th:wigner}
If $\cL$ has variance $0<\sigma^2<\infty$, then a.s.
%
%
\begin{equation}\label{eq:wigner}
\mu_{\sqrt{n}K}:=\frac{1}{n}\sum_{k=1}^n\delta_{\sqrt{n}\lambda_k(K)}
\mathop{\longrightarrow}^{w}_{n\to\infty} \mathcal{W}_{2\sigma} ,
\end{equation}
where $\mathcal{W}_{2\sigma}$ is the Wigner semi-circle law on
$[-2\sigma,+2\sigma]$. Moreover, if $\mathcal{L}$ has finite fourth moment,
then $\lambda_2(\sqrt{n}K)$ and $\lambda_n(\sqrt{n}K)$ converge
a.s. to the
edge of the limiting support $[-2\sigma,+2\sigma]$.
\end{theorem}

This Wigner-like scenario can be dramatically altered if we allow $\cL
$ to
have a heavy tail at infinity. For any $\alpha\in(0,\infty)$, we say that
$\cL$ belongs to the class $\mathbb{H}_\alpha$ if $\cL$ is
supported in
$[0,\infty)$ and has a regularly varying tail of index $\alpha$, that
is, for
all $t > 0$,
%
%
\begin{equation}\label{htp}
G(t):=\cL((t,\infty))=L(t) t^{-\alpha},
\end{equation}
where $L$ is a function with slow variation at $\infty$; that is, for
any $x>0$,
\[
\lim_{t\to\infty}\frac{L(x t)}{L(t)} = 1.
\]
Set $a_n = \inf\{a>0 \dvtx n G(a)\leq1\}$. Then
$nG(a_n)=nL(a_n)a_n^{-\alpha}\to1$ as $n\to\infty$, and
%
%
\begin{equation}\label{ht1}
n G(a_n t)\to t^{-\alpha} %
\qquad\mbox{as } %
n\to\infty%
\mbox{ for all $t > 0$}.
\end{equation}
It is well known that $a_n$ has regular variation at $\infty$ with index
$1/\alpha$, that is,
\[
a_n=n^{1/\alpha}\ell(n)
\]
for some function $\ell$ with slow variation at $\infty$ (see, e.g.,
Resnick~\cite{resnick}, Section~2.2.1). As an example, if $V$ is uniformly
distributed on the interval $[0,1]$, then for every $\alpha\in
(0,\infty)$, the
law of $V^{-1/\alpha}$, supported in $[1,\infty)$, belongs to
$\mathbb{H}_\alpha$. In this case, $L(t)=1$ for $t\geq1$, and
$a_n=n^{1/\alpha}$.

To understand the limiting behavior of the spectrum of $K$ in the heavy-tailed
case it is important to consider first the symmetric i.i.d. matrix
corresponding to the un-normalized weights $U_{i,j}$. More generally, we
introduce the random $n\times n$ symmetric matrix $X$ defined by
%
%
\begin{equation}\label{levymatrix}
X = (X_{i,j})_{1\leq{i,j}\leq n},
\end{equation}
where $(X_{i,j})_{1\leq i\leq j\leq n}$ are i.i.d. such that
$U_{i,j}:=|X_{i,j}|$
has law in $\mathbb{H}_\alpha$ with $\alpha\in(0,2)$, and
%
%
\begin{equation}\label{theta}
\theta= \lim_{t \to\infty} %
\frac{\mathbb{P}(X_{i,j}>t)}{\mathbb{P}(|X_{i,j}|>t)} \in[0,1] .
\end{equation}
It is well known that, for $\alpha\in(0,2)$, a random variable $Y$ is
in the
domain of attraction of an $\alpha$-stable law iff the law of $|Y|$ is in
$\mathbb{H}_\alpha$ and the limit (\ref{theta}) exists (cf.~\cite{Feller},
Theorem IX.8.1a). It will be useful to view the entries $X_{i,j}$ in
(\ref{levymatrix}) as the marks across edge $\{i,j\}\in E_n$ of a random
network $(G_n,\mathbf{X})$, just as the marks $U_{i,j}$ defined the network
$(G_n,\mathbf{U})$ introduced above.

Remarkable works have been devoted recently to the asymptotic behavior
of the
ESD of matrices $X$ defined by (\ref{levymatrix}), sometimes called L\'{e}vy
matrices.
The analysis of the \textit{Limiting Spectral Distribution} (LSD) for
$\alpha\in(0,2)$ is considerably harder than the finite second moment
case 
(Wigner matrices), and the LSD is nonexplicit.
Theorem~\ref{th:iida} below has been investigated by the physicists Bouchaud
and Cizeau~\cite{BouchaudCizeau} and rigorously proved by Ben Arous and
Guionnet~\cite{benarous-guionnet}, and Belinschi, Dembo and Guionnet
\cite{belinschi} (see also Zakharevich~\cite{zakharevich} for related
results).
%
\begin{theorem}[{[Symmetric i.i.d. matrix, $\alpha\in(0,2)$]}]\label{th:iida}
For every $\alpha\in(0,2)$, there exists a symmetric probability
distribution $\mu_\alpha$ on $\mathbb{R}$ depending only on $\alpha
$ such
that [with the notation of (\ref{ht1}) and (\ref{levymatrix})] a.s.
\[
\mu_{a_n^{-1}X}:=\frac1n\sum_{i=1}^n\delta_{\lambda_i(a_n^{-1}X)}
\mathop{\longrightarrow}^{w}_{n\to\infty}
\mu_\alpha.
\]
\end{theorem}

In Section~\ref{sec:iida}, we give a new independent proof of Theorem
\ref{th:iida}. The key idea of our proof is to exhibit a limiting
self-adjoint operator $\bT$ for the sequence of matrices $a_n^{-1} X$,
defined on a suitable Hilbert space, and then use known spectral convergence
theorems of operators. The limiting operator will be defined as the
``adjacency matrix'' of an infinite rooted tree with random edge
weights, the
so-called Poisson weighted infinite tree (PWIT) introduced by Aldous
\cite{aldous92} (see also~\cite{aldoussteele}). In other words, the
PWIT will
be shown to be the local weak limit of the random network $(G_n,\mathbf
{X})$ when
the edge marks $X_{i,j}$ are rescaled by $a_n$. In this setting the LSD
$\mu_\alpha$ arises as the expected value of the (random) spectral
measure of
the operator $\bT$ at the root of the tree. The PWIT and the limiting operator
$\bT$ are defined in Section~\ref{sec:PWIT}. Our method of proof can
be seen
as a variant of the resolvent method, based on local convergence of operators.
It is also well suited to investigate properties of the LSD $\mu
_\alpha$ (cf.
Theorem~\ref{th:mua} below).

Let us now come back to our random reversible Markov kernel $K$ defined by
(\ref{srw}) from weights with law $\cL\in\mathbb{H}_\alpha$. We obtain
different limiting behavior in the two regimes $\alpha\in(0,1)$ and
$\alpha\in(1,2)$. The case $\alpha>2$ corresponds to a Wigner-type behavior
(special case of Theorem~\ref{th:wigner}). We set
\[
\kappa_n=na_n^{-1}.
\]

\begin{theorem}[{[Reversible Markov matrix, $\alpha\in(1,2)$]}]\label{th:k12}
Let $\mu_\alpha$ be the probability distribution which appears as the
LSD in
the symmetric i.i.d. case (Theorem~\ref{th:iida}). If
$\cL\in\mathbb{H}_\alpha$ with $\alpha\in(1,2)$ then a.s.
\[
\mu_{\kappa_n K} :=\frac{1}{n}\sum_{k=1}^n\delta_{\lambda
_k(\kappa_n K)}%
\mathop{\longrightarrow}^{w}_{n\to\infty}\mu_\alpha.
\]
\end{theorem}
\begin{theorem}[{[Reversible Markov matrix, $\alpha\in(0,1)$]}]\label{th:k01}
For every $\alpha\in(0,1)$, there exists a symmetric probability
distribution $\wt\mu_\alpha$ supported on $[-1,1]$ depending only on
$\alpha$ such that a.s.
\[
\mu_{K} := \frac{1}{n}\sum_{k=1}^n\delta_{\lambda_k(K)}%
\mathop{\longrightarrow}^{w}_{n\to\infty}\wt\mu_\alpha.
\]
\end{theorem}

The proofs of Theorems~\ref{th:k12} and~\ref{th:k01} are given in
Sections~\ref{sec:k12} and~\ref{sec:k01}, respectively. As in the
proof of
Theorem~\ref{th:iida}, the main idea is to exploit convergence of our matrices
to suitable operators defined on the PWIT. To understand the scaling in
Theorem~\ref{th:k12}, we recall that if $\alpha>1$, then by the
strong law of
large numbers, we have $n^{-1}\rho_i\to1$ a.s. for every row sum
$\rho_i$,
and this is shown to remove, in the limit $n\to\infty$, all
dependencies in
the matrix $na_n^{-1} K$, so that we obtain the same behavior of the
i.i.d.
matrix of Theorem~\ref{th:iida}. On the other hand, when $\alpha\in(0,1)$,
both the sum $\rho_i$ and the maximum of its elements are on scale
$a_n$. The
proof of Theorem~\ref{th:k01} shows that the matrix $K$ converges (without
rescaling) to a random stochastic self-adjoint operator $\bK$ defined
on the
PWIT. The operator $\bK$ can be described as the transition matrix of the
simple random walk on the PWIT and is naturally linked to Poisson--Dirichlet
random variables. This is based on the observation that the order statistics
of any given row of the matrix $K$ converges weakly to the Poisson--Dirichlet
law $\mathrm{PD}(\alpha,0)$ (see Lemma~\ref{le:PoiExt} below for the
details).
In fact, the operator $\bK$ provides an interesting generalization of the
Poisson--Dirichlet law.

Since $\mu_K$ is supported in $[-1,1]$, (\ref{moms}) and Theorem~\ref{th:k01}
imply that for all $\ell\geq1$, a.s.
%
%
\begin{equation}\label{momms}
\frac1n\sum_{i =1} ^n p_\ell(i)
= \int_\mathbb{R} x^\ell\mu_K(dx) %
\mathop{\longrightarrow}_{n\to\infty} %
\int_\mathbb{R} x^\ell\wt\mu_\alpha(dx)=:\gamma_\ell.
\end{equation}
The LSD $\wt\mu_\alpha$ will be obtained as the expectation of the (random)
spectral measure of $\bK$ at the root of the PWIT. It will follow that
$\gamma_\ell$ (the $\ell$th moment of $\wt\mu_\alpha$) is the\vadjust{\goodbreak}
expected value of the (random) probability that the random
walk 
returns to the root in $\ell$-steps. In particular, the symmetry of
$\wt\mu_\alpha$ follows from the bipartite nature of the PWIT.

It was proved by Ben Arous and Guionnet~\cite{benarous-guionnet},
Remark 1.5,
that $\alpha\in(0,2)\mapsto\mu_\alpha$ is continuous with respect
to weak
convergence of probability measures, and by Belinschi, Dembo and Guionnet
\cite{belinschi}, Remark 1.2 and Lemma 5.2, that $\mu_\alpha$ tends
to the
Wigner semi-circle law as $\alpha\nearrow2$. We believe that Theorem
\ref{th:k12} should remain valid for $\alpha=2$ with LSD given by the Wigner
semi-circle law. Further properties of the measures $\mu_\alpha$ and
$\wt\mu_\alpha$ are discussed below.

The case $\alpha=1$ is qualitatively similar to the case $\alpha\in
(1,2)$ with
the difference that the sequence $\kappa_n$ in Theorem~\ref{th:k12}
has to be
replaced by $\kappa_n=na_n^{-1}w_n$ where
%
%
\begin{equation}\label{wwnn}
w_n=\int_0^{a_n} x\cL(dx) .
\end{equation}
Indeed, here the mean of $U_{i,j}$ may be infinite and the closest one
gets to a
law of large numbers is the statement that $\rho_i/ nw_n \to1$ in probability
(see Section~\ref{sec:k1}). The sequence $w_n$ (and therefore $\kappa
_n$) is
known to be slowly varying at $\infty$ for $\alpha=1$ (see, e.g., Feller
\cite{Feller}, VIII.8). The following mild condition will be assumed:
%
There exists $0 < \varepsilon< 1/2$ such that
%
%
\begin{equation}\label{eq:conda1wn}
\liminf_{n\to\infty}\frac{ w_{\lfloor n^\varepsilon\rfloor} }{
w_n} > 0.
\end{equation}
For example, if $U_{i,j}^{-1}$ is uniform on $[0,1]$, then $\kappa
_n=w_n = \log
n $ and\break $\lim_{n\to\infty} w_{\lfloor n^\varepsilon\rfloor} / w_n =
\varepsilon$.
In the next theorem $\mu_1$ stands for the LSD $\mu_\alpha$ from Theorem
\ref{th:iida}, at $\alpha=1$.
%
\begin{theorem}[(Reversible Markov matrix, $\alpha=1$)]\label{th:k1}
Suppose that $\cL\in\mathbb{H}_\alpha$ with $\alpha=1$ and assume
(\ref{eq:conda1wn}). If $\mu_{\kappa_n K}$ is the ESD of $\kappa
_nK$, with
$\kappa_n=na_n^{-1}w_n$, then, as $n\to\infty$, a.s. $\mu_{\kappa
_n K}%
\mathop{\longrightarrow}\limits^{w}\limits_{n\to\infty}
\mu_1$.
\end{theorem}

\subsubsection*{Properties of the LSD}
In Section~\ref{sec:spec} we prove some properties of the LSDs $\mu
_\alpha$
and $\wt\mu_\alpha$.
\begin{theorem}[(Properties of $\mu_\alpha$)]\label{th:mua}
Let $\mu_\alpha$ be the symmetric LSD in Theorems~\ref{th:iida} and \ref
{th:k12}.
\begin{enumerate}[(iii)]
\item[(i)] $\mu_\alpha$ is absolutely continuous on $\mathbb{R}$.
\item[(ii)] The density of $\mu_\alpha$ at $0$ is equal to
\[
\frac{1}{\pi} \Gamma\biggl(1+ \frac{2}{\alpha}\biggr) \biggl(
\frac{
\Gamma(1 - {\alpha}/{2} )}{ \Gamma(1 + {\alpha}/{2} ) }
\biggr)^{{1}/{\alpha}} .
\]
\item[(iii)] $\mu_\alpha$ is heavy tailed, and as $t$ goes to
$+\infty$,
\[
\mu_\alpha(( t , + \infty)) \sim\tfrac{1}{2}t^{- \alpha}.\vadjust{\goodbreak}
\]
\end{enumerate}
\end{theorem}

Statements (i) and (ii) answer some questions raised in
\cite{benarous-guionnet,belinschi}. Statement~(iii) is already
contained in
\cite{belinschi}, Theorem 1.7, but we provide a new proof based on a Tauberian
theorem for the Cauchy--Stieltjes transform that may be of independent
interest.
\begin{theorem}[(Properties of $\wt\mu_\alpha$)]\label{th:muabis}
Let $\wt\mu_\alpha$ be the symmetric LSD in Theorem~\ref{th:k01}, with
moments $\gamma_\ell$ as in (\ref{momms}). Then the following statements
hold true.
\begin{enumerate}[(iii)]
\item[(i)]
For $\alpha\in(0,1)$, there exists $\delta>0$ such that
%
\[
\gamma_{2n}\geq\delta n^{-\alpha} \qquad\mbox{for all } n
\geq1 .
\]
Moreover, we have $\liminf_{\alpha\nearrow1}\gamma_2>0$.
%
\item[(ii)] For the topology of the weak convergence, the map $\alpha
\mapsto\widetilde\mu_\alpha$ is continuous in $(0,1)$.
\item[(iii)] For the topology of the weak convergence,
\[
\lim_{\alpha\searrow0} \widetilde\mu_\alpha%
= \frac1 4 \delta_{-1}+ \frac1 2 \delta_{0} + \frac1 4 \delta_{1}.
\]
\end{enumerate}
\end{theorem}

It is delicate to provide liable numerical simulations of the ESDs.
%
%
\begin{figure}

\includegraphics{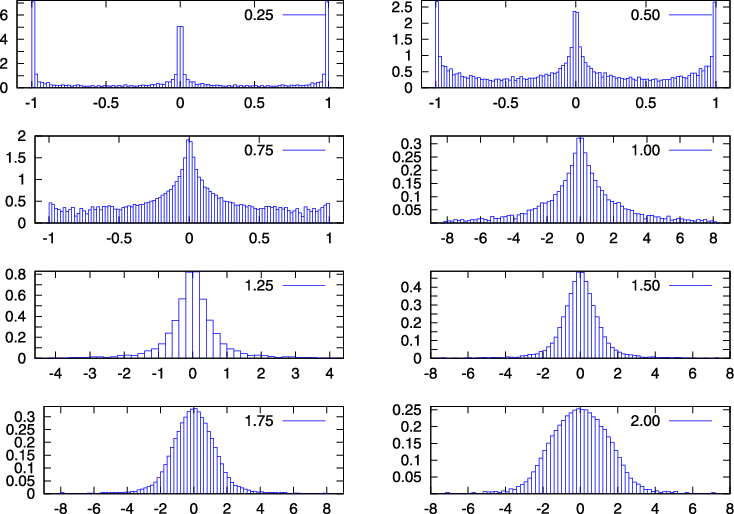}

\caption{Histograms of scaled ESDs illustrating the convergence stated by
Theorems \protect\ref{th:k12} and \protect\ref{th:k01}, for the
following values of
$\alpha\dvtx 0.25$, $0.50$, $0.75$, $1.00$, $1.25$, $1.50$, $1.75$,
$2.00$. Here $n=5000$ and $\mathcal{L}$ is the law of $V^{-1/\alpha}$
where $V$ is a uniform random variable on $(0,1)$. The first three plots
are the histogram of the spectrum of a single realization of $K$. The
fourth plot corresponds to $\alpha=1$ and is a histogram of the spectrum
of a single realization of $\log(n)K$. The four last plots are the
histogram of the spectrum of a single realization of $\kappa_n K$. In
order to avoid scaling problems, an asymptotically negligible portion of
the spectrum edge was discarded: only
$\lambda_{\lfloor\log(n)\rfloor},\ldots,\lambda_{\lfloor
n-\log(n)\rfloor}$ were used.%
}
\label{firstpic}
\end{figure}
Neverthe\-less, Figure~\ref{firstpic} provides histograms for various
values of
$\alpha$ and a large value of $n$, illustrating Theorems
\ref{th:k12}--\ref{th:muabis}.

\subsubsection*{Invariant measure and edge behavior}

Finally, we turn to the analysis of the invariant probability distribution
$\hat\rho$ for the random walk on $(G,\bU)$. This is obtained by normalizing
the vector of row sums $\rho$
\[
\hat\rho= (\rho_1 +\cdots+ \rho_n) ^{ -1} (\rho_1,\ldots, \rho
_n ).
\]
Following~\cite{bordenave-caputo-chafai}, Lemma 2.2, if $\alpha>2$, then
$n\max_{1\leq i\leq n}|\hat\rho_i - n^{-1}| \to0$ as $n\to\infty$
a.s. This
uniform strong law of large numbers does not hold in the heavy-tailed case
$\alpha\in(0,2)$: the large $n$ behavior of $\hat\rho$ is then
dictated by
the largest weights in the system.

Below we use the notation $\wt\rho=(\wt\rho_1,\ldots,\wt\rho_n)$
for the
ranked values of $\hat\rho_1,\ldots,\hat\rho_n$, so that $\wt\rho
_1\geq
\wt\rho_2\geq\cdots$ and their sum is $1$. The symbol
$\stackrel{d}{\longrightarrow}$ denotes convergence in
distribution. We refer to Section~\ref{order} for more details on weak
convergence in the space of ranked sequences and for the definition of the
Poisson--Dirichlet law $\mathrm{PD}(\alpha,0)$.
\begin{theorem}[(Invariant probability measure)]
\label{th:inv} Suppose that $\cL\in\mathbb{H}_\alpha$.
\begin{enumerate}[(ii)]
\item[(i)] If $\alpha\in(0,1)$, then
%
%
\begin{equation}\label{inv02} \wt\rho
\mathop{\longrightarrow}^{d}_{n\to\infty}
\tfrac12 (V_1,V_1,V_2,V_2,\ldots) ,
\end{equation}
where $V_1>V_2>\cdots$ stands for a Poisson--Dirichlet
$\mathrm{PD}(\alpha,0)$ random vector.
\item[(ii)] If $\alpha\in(1,2)$, then
%
%
\begin{equation}\label{inv01}
\kappa_{n(n+1)/2} \wt\rho
\mathop{\longrightarrow}^{d}_{n\to\infty}
\tfrac12(x_1,x_1,x_2,x_2,\ldots) ,
\end{equation}
where\vspace*{1pt} $x_1>x_2>\cdots$ denote the ranked points of the Poisson point
process on $(0,\infty)$ with intensity measure $\alpha x^{-\alpha-1}\,dx$.
Moreover, the same convergence holds for $\alpha=1$ provided the sequence
$\kappa_n$ is replaced by $n a_n^{-1} w_n$, with $w_n$ as in (\ref{wwnn}).
\end{enumerate}
\end{theorem}

Theorem~\ref{th:inv} is proved in Section~\ref{sec:inv}. These
results will be
derived from the statistics of the ranked values of the weights $U_{i,j}$,
$i<j$, on the scale $a_{n(n+1)/2}$ (diagonal weights $U_{i,i}$ are
easily seen
to give negligible contributions). The duplication in the sequences in
(\ref{inv01}) and (\ref{inv02}) then comes from the fact that each of the
largest weights belongs to two distinct rows and determines alone the limiting
value of the associated row sum.

Theorem~\ref{th:inv} is another indication that the random walk with
transition matrix $K$ shares the features of a \textit{trap model}. Loosely
speaking, instead of being trapped at a vertex, as in the usual mean field
trap models 
(see~\cite{Bouchaud,BenArousCerny,MR2152251,MR2435851}) here the
walker is
trapped at an edge.

Large edge weights are responsible for the large eigenvalues of $K$. This
phenomenon is well understood in the case of symmetric random matrices with
i.i.d. entries, where it is known that, for $\alpha\in(0,4)$, the
edge of the
spectrum gives rise to a Poisson statistics (see
\cite{MR2081462,auffinger-benarous-peche}). The behavior of the extremal
eigenvalues of $K$ when $\mathcal{L}$ has finite fourth moment has been
studied in~\cite{bordenave-caputo-chafai}. In particular, it is shown there
that the spectral gap $1-\lambda_2$ is $1-O(n^{-1/2})$. In the present
case of
heavy-tailed weights, in contrast, by localization on the largest
edge weight it is possible to prove that, a.s. and up to corrections with
slow variation at $\infty$,
%
%
\begin{equation}\label{edges} 1-\lambda_2 =
\cases{
O(n^{-1/\alpha}), &\quad $\alpha\in(0,1)$,\cr
O\bigl(n^{-(2-\alpha)/\alpha}\bigr), &\quad $\alpha\in[1,2)$.}
\end{equation}
Similarly, for $\alpha\in(2,4)$ one has that $\lambda_2$ is bounded
below by
$n^{-(\alpha-2)/\alpha}$. Understanding the statistics of the extremal
eigenvalues remains an interesting open problem.

\section{Convergence to the Poisson weighted infinite tree}
\label{sec:PWIT}

The aim of this section is to prove that the matrices $X$ and $K$
appearing in
Theorems~\ref{th:iida},~\ref{th:k12} and~\ref{th:k01}, when properly
rescaled,
converge ``locally'' to a limiting operator defined on the Poisson weighted
infinite tree (PWIT). The concept of local convergence of operators is defined
below. We first recall the standard construction of the PWIT.

\subsection{The PWIT}

Given a 
Radon measure $\nu$ on $\mathbb{R}$, $\pwit(\nu)$ is the random
rooted tree
defined as follows. The vertex set of the tree is identified with
$\mathbb{N}^f:= \bigcup_{k \in\mathbb{N}} \mathbb{N}^k$ by indexing
the root as
$\mathbb{N}^0 = \eset$, the offsprings of the root as $\mathbb{N}$
and, more
generally, the offsprings of some $\bv\in\mathbb{N}^k$ as $(\bv
1),(\bv2),
\ldots\in\mathbb{N}^{k+1}$ [for short notation, we write $(\bv1)$
in place
of $(\bv,1)$]. In this way the set of $\bv\in\mathbb{N}^n$
identifies the
$n$th generation.

We now assign marks to the edges of the tree according to a collection
$\{
\Xi_\bv\}_{\bv\in\mathbb{N}^f}$ of independent realizations of the Poisson
point process with intensity measure $\nu$ on $\mathbb{R}$. Namely, starting
from the root $\eset$, let $ \Xi_\eset= \{y_1,y_2,\ldots\}$ be
ordered in such
a way that $|y_1| \leq|y_2| \leq\cdots,$ and assign the mark $y_i$
to the
offspring of the root labeled $i$. Now, recursively, at each vertex
$\bv$ of
generation $k$, assign the mark $y_{\bv i}$ to the offspring labeled
$\bv i$,
where $\Xi_{\bv}=\{ y_{\bv1}, y_{\bv2}, \ldots\}$ satisfy $|y_{\bv
1}| \leq
|y_{\bv2}| \leq\cdots.$


\subsection{Local operator convergence}

We give a general formulation and later specialize to our setting. Let
$V$ be
a countable set,\vadjust{\goodbreak} and let $L^2(V)$ denote the Hilbert space defined by the
scalar product
\[
\langle\phi,\psi\rangle:= \sum_{u\in V} \bar\phi_u\psi_u
,\qquad \phi_u =
\langle\delta_u,\phi\rangle,
\]
where $\phi,\psi\in\mathbb{C}^V$ and $\delta_{u}$ denote the unit vector
with support $u$. Let $\cD$ denote the dense subset of $L^2 (V)$ of vectors
with finite support.
\begin{defi}[(Local convergence)]\label{def:convloc}
Suppose $\bS_n$ is a sequence of bound\-ed operators on $L^2(V)$, and
$\bS$ is
a closed linear operator on $L^2(V)$ with dense domain $D(\bS)\supset
\cD$.
Suppose further that $\cD$ is a core for $\bS$ (i.e., the closure of
$\bS$
restricted to $\cD$ equals $\bS$). For any $u,v\in V$ we say that
$(\bS_n,u)$ converges locally to $(\bS,v)$ and write
\[
(\bS_n,u) \to(\bS,v) ,
\]
if there exists a sequence of bijections $\sigma_n\dvtx V\to V$ such that
$\sigma_n (v) = u$ and, for all $\phi\in\cD$,
\[
\sigma_n ^{-1} \bS_n \sigma_n \phi\to\bS\phi,
\]
in $L^2(V)$, as $n\to\infty$.
\end{defi}

In other words, this is the standard strong convergence of operators up
to a
re-indexing of $V$ which preserves a distinguished element. With a slight
abuse of notation we have used the same symbol $\sigma_n$ for the linear
isometry $\si_n\dvtx L^2(V)\to L^2(V)$ induced in the obvious way, that
is, such that
$\si_n\delta_v = \delta_{\si_n(v)}$ for all $v\in V$. The point for
introducing Definition~\ref{def:convloc} lies in the following theorem on
strong resolvent convergence. Recall that if $\mathbf{S}$ is a self-adjoint
operator its spectrum is real, and for all $z \in\mathbb{C}_+ := \{z
\in
\mathbb{C} \dvtx\Im z > 0\} $, the operator $\bS- zI$ is invertible
with bounded
inverse. The operator-valued function $z \mapsto(\bS- zI)^{-1}$ is the
resolvent of~$\bS$.
\begin{theorem}[(From local convergence to resolvents)]\label{th:strongres}
If $\bS_n$ and $\bS$ are self-adjoint operators that satisfy the conditions
of Definition~\ref{def:convloc} and $(\bS_n,\break u) \to(\bS,v)$ for some
$u,v\in
V$, then, for all $z \in\mathbb{C}_+$,
%
%
\begin{equation}\label{strresconv} \langle
\delta_{u}, (\bS_n- zI) ^{-1} \delta_u \rangle\to\langle\delta_{v},
(\bS- zI) ^{-1} \delta_{v} \rangle.
\end{equation}
\end{theorem}
\begin{pf}
It is a special case of 
\cite{reedsimon}, Theorem VIII.25(a). Indeed, if we define $\wt\bS_n =
\sigma_n ^{-1} \bS_n \sigma_n$, then\vspace*{1pt} $\wt\bS_n \phi\to\bS\phi$
for all
$\phi$ in a common core of the self-adjoint operators $\wt\bS_n, \bS$.
This implies\vspace*{1pt} the strong resolvent convergence, that is, $(\wt\bS_n -
zI)^{-1}\psi\to(\bS- zI)^{-1}\psi$ for any $z\in\mathbb{C}_+$,
$\psi\in
L^2(V)$. The conclusion follows by taking the scalar product
\[
\langle\delta_v , (\wt\bS_n - zI)^{-1}\delta_v\rangle%
= \langle\delta_u , (\bS_n - zI)^{-1}\delta_u\rangle.
\]
\upqed\end{pf}

We shall apply the above theorem in cases where the operators $\bS_n$ and
$\bS$ are random operators on $L^2(V)$,\vadjust{\goodbreak} which satisfy with probability
one the
conditions of Definition~\ref{def:convloc}. In this case we say that
$(\bS_n
,u) \to(\bS,v)$ \textit{in distribution} if there exists a random bijection
$\sigma_n$ as in Definition~\ref{def:convloc} such that $\sigma_n
^{-1} \bS_n
\sigma_n \phi$ converges in distribution to $\bS\phi$, for all
$\phi\in
\cD$ [where a random vector $\psi_n \in L^2 (V)$ converges in
distribution to
$\psi$ if
\[
\lim_{n\to\infty} \mathbb{E} f (\psi_n) = \mathbb{E} f(\psi)
\]
for all bounded continuous functions $f\dvtx L^2 (V)\to\mathbb{R}$].
Under these
assumptions then (\ref{strresconv}) becomes convergence in
distribution of
(bounded) complex random variables. In our setting the Hilbert space
will be
$L^2(V)$, with $V=\mathbb{N}^f$, the vertex set of the PWIT, the operator
$\bS_n$ will be a rescaled version of the matrix $X$ defined by
(\ref{levymatrix}) or the matrix $K$ defined by (\ref{srw}). The operator
$\bS$ will be the corresponding limiting operator defined below.

\subsection{Limiting operators}\label{limop}

Let $\theta$ be as in Theorem~\ref{th:iida}, and let $\ell_\theta$
be the
positive Borel measure\vspace*{1pt} on the real line defined by $d\ell_\theta
(x)=\theta
\ind_{\{x > 0\}}\,dx + (1-\theta)\ind_{\{x < 0\}}\,dx$. Consider a
realization of
$\pwit(\ell_\theta)$. As before the mark from vertex $\bv\in
\mathbb{N}^k$ to
$\bv k \in\mathbb{N}^{k+1}$ is denoted by $y_{\bv k}$. We note that almost
surely
%
%
\begin{equation}\label{l20} \sum_k |y_{\bv k }|^{-2/\alpha} < \infty,
\end{equation}
since a.s. $\lim_k |y_{\bv k } | / k = 1$ and $\sum_k k^{-2/\alpha}$
converges for $\alpha\in(0,2)$. Recall that for $V = \mathbb{N} ^f$,
$\cD$
is the dense set of $L^2 (V)$ of vectors with finite support. We may
a.s.
define a linear operator $\bT\dvtx\cD\to L^2(V)$ by letting, for $\bv
,\bw\in
\mathbb{N}^f$,
%
%
\begin{eqnarray}\label{tone}\qquad
\bT(\bv,\bw) &=& \langle\delta_{\bv} , \bT\delta_{\bw} \rangle\nonumber\\[-8pt]\\[-8pt]
&=&
\cases{
\sign(y_{\bw}) |y_{\bw}|^{-1/\alpha},
&\quad if $\bw= \bv k$ for some integer $k$, \cr
\sign(y_{\bv}) |y_{\bv}|^{-1/\alpha},
&\quad if $\bv= \bw k$ for some integer $k$,\cr
0, &\quad otherwise.}\nonumber
\end{eqnarray}
Note that if every edge $e$ in the tree with mark $y_e$ is given the
``weight'' $\sign(y_{e}) |y_{e}|^{-1/\alpha}$ then we may look at the operator
$\bT$ as the ``adjacency matrix'' of the weighted tree. Clearly, $\bT
$ is
symmetric, and therefore it has a closed extension with domain $D(\bT)
\subset
L^2 (\mathbb{N}^f)$ such that $\cD\subset D(\bT)$ (see,
e.g., 
\cite{reedsimon}, Chapter~VIII, Section 2). We will prove in Proposition~\ref{esssa}
below that
$\bT$ is essentially self-adjoint, that is, the closure of $\bT$ is
self-adjoint. With a slight abuse of notation, we identify $\bT$ with its
closed extension. As stated below, $\bT$ is the weak local limit of the
sequence of $n\times n$ i.i.d. matrices $a_n^{-1}X$, where $X$ is
defined by
(\ref{levymatrix}). To this end we view the matrix $X$ as an operator in
$L^2(V)$ by setting $\langle\delta_i , X\delta_j\rangle= X_{i,j}$, where
$i,j\in\mathbb{N}$ denote the labels of the offsprings of the root
(the first
generation), with the convention that $X_{i,j}=0$ when either $i>n$ or $j>n$,
and by setting $\langle\delta_\mathbf{u} , X\delta_\bv\rangle=
0$ when
either $\mathbf{u}$ or
$\bv$ does not belong to the first generation.\vadjust{\goodbreak}

Similarly, taking now $\theta=1$, in the case of Markov matrices $K$ defined
by (\ref{srw}), for $\alpha\in[1,2)$, $\bT$ is the local limit
operator of
$\kappa_n K$. To work directly with symmetric operators we introduce the
symmetric matrix
%
%
\begin{equation}\label{ksim}
S_{i,j} = \frac{U_{i,j}}{\sqrt{\rho_i\rho_j}} ,
\end{equation}
which is easily seen to have the same spectrum of $K$ (see, e.g.,
\cite{bordenave-caputo-chafai}, Lemma~2.1). Again the matrix $S$ can be embedded in the
infinite tree as described above for $X$.

In the case $\alpha\in(0,1)$ the Markov matrix $K$ has a different limiting
object that is defined as follows. Consider a realization of $\pwit
(\ell_1)$,
where $\ell_1$ is the Lebesgue measure on $[0,\infty)$. We define an operator
corresponding to the random walk on this tree with conductance equal to the
mark to the power $-1/\alpha$. More precisely, for $\bv\in\mathbb
{N}^f$, let
\[
\rho(\bv) = y^{-1/\alpha}_{\bv} + \sum_{k \in\mathbb{N}}
y^{-1/\alpha}_{\bv k}
\]
with the convention that $y^{-1/\alpha}_{\eset} = 0$. Since a.s.
$\lim_k |y_{\bv
k } | / k = 1$, $\rho(\bv)$ is almost surely finite for $\alpha\in
(0,1)$. We
define the linear operator $\bK$ on $\cD$, by letting, for $\bv,\bw
\in
\mathbb{N}^f$,
%
%
\begin{equation}\label{kappone}\qquad
\bK(\bv,\bw) %
= \langle\delta_{\bv} , \bK\delta_{\bw} \rangle%
= \cases{
\dfrac{y^{-1/\alpha}_{\bw}}{\rho(\bv)},
&\quad if $\bw= \bv k$ for some integer $k$,\vspace*{2pt}\cr
\dfrac{y^{-1/\alpha}_{\bv}}{\rho(\bv)},
&\quad if $\bv= \bw k$ for some integer $k$,\vspace*{2pt}\cr
0, &\quad otherwise.}
\end{equation}
Note that $\bK$ is not symmetric, but it becomes symmetric in the weighted
Hilbert space $L^2(V,\rho)$ defined by the scalar product
\[
\langle\phi, \psi\rangle_\rho:=
\sum_{\mathbf{u}\in V} \rho(\mathbf{u})
\bar\phi_\mathbf{u}\psi_\mathbf{u} .
\]
Moreover, on $L^2(V,\rho)$, $\bK$ is a bounded self-adjoint operator since
Schwarz's inequality implies
\begin{eqnarray*}
\langle\bK\phi, \bK\phi\rangle_\rho^2 &=&
\sum_\mathbf{u}\rho(\mathbf{u})
\biggl|\sum_\bv\bK(\mathbf{u},\bv)
\phi_\bv\biggr|^2\\
&\leq&
\sum_\mathbf{u}\rho(\mathbf{u})\sum_\bv\bK(\mathbf{u},\bv
)|\phi_\bv|^2\\
&=& \sum_\bv\rho(\bv)
|\phi_\bv|^2
=\langle\phi, \phi\rangle_{\rho}^2
\end{eqnarray*}
so that the operator norm of $\bK$ is less than or equal to $1$. To
work with
self-adjoint operators in the unweighted Hilbert space $L^2(V)$ we shall
actually consider the operator $\bS$ defined by
%
%
\begin{equation}\label{opsym}
\bS(\bv,\bw):=\sqrt{\frac{\rho(\bv)}{\rho(\bw)}} \bK(\bv
,\bw) =
\frac{\bT(\bv,\bw)} {\sqrt{\rho(\bv) \rho(\bw)}} .
\end{equation}
This defines a bounded self-adjoint operator in $L^2(V)$. Indeed, the map
$\delta_\bv\to\sqrt{\rho(\bv)}\delta_\bv$ induces a linear
isometry $\bD\dvtx
L^2(V,\rho)\to L^2(V)$ such that
%
%
\begin{equation}\label{isom}
\langle\phi, \bS\psi\rangle= \langle\bD^{-1}\phi, \bK\bD
^{-1}\psi\rangle_\rho,
\end{equation}
for all $\phi,\psi\in L^2(V)$. In this way, when $\alpha\in(0,1)$,
$\bS$ will be
the limiting operator associated with the matrix $S$ defined in (\ref{ksim}).
Note that no rescaling is needed here. The main result of this section
is the
following.
\begin{theorem}[(Limiting operators)]\label{th:convop}
As $n$ goes to infinity, in distribution:
\begin{enumerate}[(iii)]
\item[(i)] if $\alpha\in(0,2)$ and $\theta\in[0,1]$, then $(a_n
^{-1} X
,1) \to(\bT,\eset)$;
\item[(ii)] if $\alpha\in(1,2)$ and $\theta= 1$, then $( \kappa_n
S ,1)
\to(\bT,\eset)$;
\item[(iii)] if $\alpha\in(0,1)$, then $(S ,1) \to(\bS,\eset)$.
\end{enumerate}
\end{theorem}

From the remark after Theorem~\ref{th:strongres} we see that Theorem
\ref{th:convop} implies convergence in distribution of the resolvent
at the
root. As we shall see in Section~\ref{se:convesd}, this in turn gives
convergence of the expected values of the Cauchy--Stieltjes transform
of the
ESD of our matrices.
The rest of this section is devoted to the proof of Theorem~\ref{th:convop}.

\subsection{Weak convergence of a single row}\label{order}

In this paragraph, we recall some facts about the order statistics of the
first row of the matrix $X$ and $K$, that~is,
%
\[
(X_{1,1},\ldots,X_{1,n}) \quad\mbox{and}\quad 
(U_{1,1},\ldots,U_{1,n})/\rho_1 ,
\]
where $U_{1,j} = |X_{1,j}|$ has law $\mathbb{H}_\alpha$.
Let us
denote by $V_1\geq V_2\geq\cdots\geq V_n$ the order statistics of the
variables $U_{1,j}$, $1\leq j\leq n$. Recall that $\rho_1 =
\sum_{j=1}^nV_j$. Let us define $\Delta_{k,n}=\sum_{j=k+1}^nV_j$ for
$k<n$ and
$\Delta^2_{k,n}=\sum_{j=k+1}^nV_j^2$. Call $\cA$ the set of
sequences $\{v_j\}\in
[0,\infty)^\mathbb{N}$ with $v_1\geq v_2\geq\cdots\geq0$ such that
$\lim_{j\to\infty} v_j=0$, and let $\cA_1\subset\cA$ be the
subset of
sequences satisfying $ \sum_j v_j = 1$. We shall view
\[
Y_n = \biggl(\frac{V_1}{a_n},\ldots,\frac{V_n}{a_n}\biggr) %
\quad\mbox{and}\quad %
Z_n = \biggl(\frac{V_1}{\rho_1},\ldots,\frac{V_n}{\rho_1}\biggr)
\]
as elements of $\cA$ and $\cA_1$, respectively, simply by adding
zeros to the
right of $V_n/a_n$ and $V_n/\rho_1$. Equipped with the standard product
metric, $\cA$ and $\cA_1$ are complete separable metric spaces ($\cA
_1$ is
compact), and convergence in distribution for $\cA,\cA_1$-valued random
variables is equivalent to finite-dimensional convergence (cf., e.g., Bertoin
\cite{bertoin06}).\vadjust{\goodbreak}

Let $E_1,E_2,\ldots$ denote i.i.d. exponential variables with mean
$1$ and write $\gamma_k=\sum_{j=1}^k E_j$. We define the random
variable in $\cA$
\[
Y= (\gamma_1^{-1/\alpha},\gamma_2^{-1/\alpha},
\ldots).
\]
The law of $Y$ is the law of the ordered points of a Poisson process on
$(0,\infty)$ with intensity measure $\alpha x ^{-\alpha- 1} \,dx$. For
$\alpha
\in(0,1)$ we define the variable in~$\cA_1$
\[
Z=\biggl(\frac{\gamma_1^{-1/\alpha}}{\sum_{n=1}^\infty\gamma
_n^{-1/\alpha}} , \frac{\gamma_2^{-1/\alpha}}{\sum
_{n=1}^\infty\gamma_n^{-1/\alpha}} , \ldots\biggr) .
\]
For $\alpha\in(0,1)$ the sum $\sum_n\gamma_n^{-1/\alpha}$ is
a.s. finite.
The law of $Z$ in $\cA_1$ is called the \textit{Poisson--Dirichlet} law
$\mathrm{PD}(\alpha,0)$ (see Pitman and Yor~\cite{MR1434129}, Proposition
10).
The next result is rather standard but we give a simple proof for convenience.

\begin{lem}[(Poisson--Dirichlet laws and Poisson point processes)]%
\label{le:PoiExt}
\begin{enumerate}[(iii)]
\item[(i)] For all $\alpha>0$, $Y_n$ converges in distribution to $Y$.
Moreover, for $\alpha\in(0,2)$, $(a_n^{-1} V_j)_{j \geq1}$ is a.s.
uniformly square integrable, that is, a.s.\break $\lim_{k} \sup_{n>
k} a_n^{-2}\times \Delta^2_{k,n} = 0$.
\item[(ii)] If $\alpha\in(0,1)$, $Z_n$ converges in distribution to $Z$.
Moreover, $(a_n^{-1} V_j)_{j \geq1}$ is a.s. uniformly integrable,
that is,
a.s. $\lim_{k} \sup_{n> k} a_n^{-1} \Delta_{k,n} = 0$.
\item[(iii)] If $I \subset\mathbb{N}$ is a finite set and
$V^I_{1}\geq
V^I_{2}\geq\cdots$ denote the order statistics of $\{U_{1,j}\}_{j \in
\{
1 ,\ldots,n \} \setminus I}$ then \textup{(i)} and \textup{(ii)} hold with $Y^I_n =
(V^I_1/a_n,V^I_2/a_n, \ldots)$ and $Z^I_n = (V^I_1/\rho_1,V^I_2/\rho_1,
\ldots)$.
\end{enumerate}
\end{lem}

As an example, from (i), we retrieve the well-known fact that for any
$\alpha>0$, the random variable $ a_n^{-1}\max(U_{1,1},\ldots
,U_{1,n}) $
converges weakly as $n\to\infty$ to the law of $\gamma_1^{-1/\alpha}$.
This law, known as a Fr\'{e}chet law, has density $\alpha
x^{-\alpha-1}e^{-x^{-\alpha}}$ on $(0,\infty)$.
\begin{pf*}{Proof of Lemma~\ref{le:PoiExt}}
As in LePage, Woodroofe and Zinn~\cite{zinn81} we take advantage of the
following well-known representation for the order statistics of i.i.d.
random variables. Let $G$ be the function in (\ref{htp}) and write
\[
G^{-1}(u) = \inf\{y>0\dvtx G(y)\leq u\} ,
\]
$u\in(0,1)$. We have that $(V_1,\ldots, V_n)$ equals in distribution the
vector
%
%
\begin{equation}\label{rep}
\bigl(G^{-1}(\gamma_1/\gamma_{n+1}),\ldots,G^{-1}
(\gamma_n/\gamma_{n+1})\bigr) ,
\end{equation}
where $\gamma_j$ has been defined above. To prove (i) we start from the
distributional identity
\[
Y_n\stackrel{d}{=}
\biggl(\frac{G^{-1}(\gamma_1/\gamma_{n+1})}{a_n},\ldots
,\frac{G^{-1}(\gamma_n/\gamma_{n+1})}{a_n}\biggr) ,
\]
which follows from (\ref{rep}). It suffices to prove that for every $k$,
almost surely the first $k$ terms above converge to the first $k$ terms in
$Y$. Thanks to (\ref{ht1}), almost surely, for every $j$,
%
%
\begin{equation}\label{conn}
a_n^{-1} G^{-1}(\gamma_j/\gamma_{n+1})\to\gamma
_j^{-1/\alpha} ,
\end{equation}
and the convergence in distribution of $Y_n$ to $Y$ follows. Moreover, from
(\ref{ht1}), for any $\delta>0$ we can find $n_0$ such that
\[
a_n^{-1} V_j = a_n^{-1} G^{-1}(\gamma_j/\gamma_{n+1}
)\leq
\bigl(n\gamma_j/(1+\delta)\gamma_{n+1}\bigr)^{-1/\alpha} ,
\]
for\vspace*{1pt} $n\geq n_0$, $j\in\mathbb{N}$. Since $n/\gamma_{n+1}\to1$,
a.s. we see
that the expression above is a.s. bounded by\vspace*{-1pt}
$2(1+\delta)^{1/\alpha}\gamma_j^{-1/\alpha}$, for $n$ sufficiently
large, and the second part of (i) follows from a.s. summability of
$\gamma_j^{-2/\alpha}$.

Similarly, if $\alpha\in(0,1)$, $\Delta_{k,n}$ has the same law of
\[
\sum_{j
=k+1}^n G^{-1}(\gamma_j/\gamma_{n+1}),
\]
and the second part of (ii) follows from a.s. summability of
$\gamma_j^{-1/\alpha}$. To prove the convergence of $Z_n$ we use the
distributional identity
\[
Z_n\stackrel{d}{=}
\biggl(\frac{G^{-1}(\gamma_1/\gamma_{n+1})}{\sum_{j=1}^n
G^{-1}(\gamma_j/\gamma_{n+1})},\ldots,\frac{G^{-1}
(\gamma_n/\gamma_{n+1})}{\sum_{j=1}^n
G^{-1}(\gamma_j/\gamma_{n+1})}\biggr) .
\]
As a consequence of (\ref{conn}), we then have almost surely
\[
a_n^{-1} \sum_{j=1}^nG^{-1}(\gamma_j/\gamma_{n+1})\to
\sum_{j=1}^\infty\gamma_j^{-1/\alpha} ,
\]
and (ii) follows. Finally, (iii) is an easy consequence of the
exchangeability of the variable $(U_{1,i})$
\[
\mathbb{P}( V^I _k \neq V_k ) %
\leq\mathbb{P} ( \exists j \in I \dvtx U_{1,j} \geq V_k ) %
\leq|I| \mathbb{P} (U_{1,1} \geq V_k) = |I|\frac{k}{n}.
\]
\upqed\end{pf*}

The intensity measure $\alpha x ^{-\alpha- 1} \,dx$ on $(0,\infty)$ is not
locally finite at $0$. It will be more convenient to work with Radon
(i.e., locally finite) intensity measures.
\begin{lem}[(Poisson point processes with Radon intensity
measures)]\label{le:HT}
Let $\xi^{n}_1,\xi^{n}_2,\ldots$ be sequences of i.i.d. random
variables on
$\overline{\mathbb{R}}:=\mathbb{R}\cup\{\pm\infty\}$ such that
%
%
\begin{equation}\label{asso1}
n \mathbb{P} ( \xi^{n}_1 \in\cdot) 
\mathop{\longrightarrow}^{w}_{n\to\infty}
\nu,
\end{equation}
where $\nu$ is a Radon measure on $\mathbb{R}$. Then, for any finite
set $I
\subset\mathbb{N} $ the random measure
\[
\sum_{i\in\{1,\ldots,n\}\setminus I} \delta_{\xi^{n}_i}
\]
converges weakly as $n\to\infty$ to $\operatorname{PPP}(\nu)$, the Poisson point
process on $\mathbb{R}$ with intensity law $\nu$, for the usual vague
topology on Radon measures.
\end{lem}

We refer to~\cite{resnick}, Theorem 5.3, page 138, for a proof of Lemma
\ref{le:HT}. Note that for $\xi^{(n)}_j = a_n / U_{1,j}$ it is a consequence
of Lemma~\ref{le:PoiExt}(iii). In the case $\xi^{(n)}_j = a_n / X_{1,j}$,
where $X_{i,j}$ is as in (\ref{levymatrix}) and (\ref{theta}), the
above lemma
yields convergence to PPP($\nu_{\alpha,\theta}$), where
%
%
\begin{equation}\label{appl}
\nu_{\alpha,\theta}(dx) = \bigl[\theta\ind_{\{x > 0\}}+(1-\theta
)\ind_{\{x < 0\}}\bigr] \alpha|x|^{\alpha-1} \,dx .
\end{equation}

\subsection{Local weak convergence to PWIT}
\label{subsec:LWC}

In the previous paragraph we have considered the convergence of the
first row
of the matrix $a_n^{-1}X$.
Here we generalize this by characterizing the limiting local
structure of the complete graph with marks $a_n / X_{i,j}$. Our
argument is
based on a technical generalization of an argument borrowed from Aldous
\cite{aldous92}. This will lead us to Theorems~\ref{th:convop} and
\ref{th:convop2} below.

Let $G_n$ be the complete network on $\{1,\ldots,n\}$ whose mark on
edge $(i,j)$ equals $\xi^n_{i,j}$, for some collection $(\xi^n_{i
j})_{1 \leq i \leq j \leq n}$ of i.i.d. random variables with values in
$\mathbb {R}$, with $\xi^n _ {j,i} = \xi^n_{i,j}$. We
consider\vspace*{1pt} the rooted network $(G_n,1)$ obtained by
distinguishing the vertex labeled $1$.

We follow Aldous~\cite{aldous92}, Section 3. For every fixed
realization of
the marks $(\xi^n_{i j})$, and for any $B,H\in\mathbb{N}$, such that
$(B^{H+1} -
1)/(B-1) \leq n$, we define a finite rooted subnetwork $(G_n,1)^{B,H}$ of
$(G_n,1)$, whose vertex set coincides with a $B$-ary tree of depth $H$ with
root at $1$.

To this end we partially index the vertices of $(G_n,1)$ as elements in
\[
J_{B,H} = \bigcup_{\ell=0}^H \{1,\ldots, B\}^\ell\subset\mathbb{N}^f,
\]
the indexing being given by an injective map $\sigma_n$ from $J_{B,H}$ to
$V_n:=\{1,\ldots,n\}$. The map $\si_n$ can be extended to a bijection
from a
subset of $\mathbb{N}^f$ to $V_n$. We set $I_\eset= \{ 1 \}$ and the
index of
the root $1$ is $ \sigma_n^{-1} (1) = \eset$. The vertex $v\in V_n
\setminus
I_{\eset}$ is given the index $(k) = \sigma_n^{-1} (v)$, $1 \leq
k\leq B$, if
$\xi^n_{(1,v)}$ has the $k$th smallest absolute value among
$\{\xi^n_{1,j}, j\neq1\}$, the marks of edges emanating from the
root $1$.
We break ties by using the lexicographic order. This defines the first
generation. Now let $I_1$ be the union of $I_\eset$ and the $B$
vertices that
have been selected. If $H\geq2$, we repeat the indexing procedure for the
vertex indexed by $(1)$ (the first child)\vadjust{\goodbreak} on the set $V_n \setminus
I_1$. We
obtain a new set $\{11,\ldots,1B\}$ of vertices sorted by their
weights as
before [for short notation, we concatenate the vector $(1,1)$ into
$11$]. Then
we define $I_{2}$ as the union of $I_1$ and this new collection. We
repeat the
procedure for $(2)$ on $V_n \setminus I_{2}$ and obtain a new set
$\{21,\ldots,2 B\}$, and so on. When we have constructed $\{B1,\ldots
,BB\}$,
we have finished the second generation (depth $2$) and we have indexed $(B^{3}
- 1)/(B-1)$ vertices. The indexing procedure is then repeated until
depth $H$
so that $(B^{H+1} - 1)/(B-1)$ vertices are sorted. Call this set of vertices
$V_n^{B,H} = \sigma_n J_{B,H} $. The\vspace*{1pt} subnetwork of $G_n$ generated by
$V_n^{B,H}$ is denoted $(G_n,1)^{B,H}$ (it can be identified with the original
network $G_n$ where any edge $e$ touching the complement of $V_n^{B,H}$ is
given a mark $x_e=\infty$). In $(G_n,1)^{B,H}$, the set $\{\mathbf
{u}1,\ldots,\mathbf{u}B
\}$ is called the set of children or offsprings of the vertex $\mathbf
{u}$. Note that
while the vertex set has been given a tree structure, $(G_n,1)^{B,H}$
is still
a complete network. The next proposition shows that it nevertheless converges
to a tree (i.e., all circuits vanish, or equivalently, the extra marks diverge
to $\infty$) if the $\xi^n_{i,j}$ satisfy a suitable scaling assumption.

Let $(\cT,\eset)$ denote the infinite random rooted network with
distribution
$\pwit(\nu)$. We call $(\cT,\eset)^{B,H}$ the finite random network obtained
by the sorting procedure described in the previous paragraph. Namely,
$(\cT,\eset)^{B,H}$ consists of the sub-tree with vertices of the
form $\mathbf{u}\in
J_{B,H}$, with the marks inherited from the infinite tree. If an edge
is not
present in $(\cT,\eset)^{B,H}$, we assign to it the mark $+\infty$.

We say that the sequence of random finite networks $(G_n,1)^{B,H}$ converges
in distribution (as $n\to\infty$) to the random finite network
$(\cT,\eset)^{B,H}$ if the joint distributions of the marks converge
weakly. To make this precise we have to add the points $\{\pm\infty\}
$ as
possible values for each mark, and continuous functions on the space of marks
have to be understood as functions such that the limit as any one of
the marks
diverges to $+\infty$ exists and coincides with the limit as the same mark
diverges to $-\infty$. The next proposition generalizes \cite
{aldous92}, Section 3.
\begin{prop}[(Local weak convergence to a tree)] \label{prop:LWC} %
Let $(\xi^n_{i,j})_{1 \leq i \leq j \leq n}$ be a collection of
i.i.d. random variables with values in
$\overline{\mathbb{R}}:=\mathbb{R}\cup\{\pm\infty\}$ and set $\xi
^n _ {j,i}
= \xi^n_{i,j}$. Let $\nu$ be a Radon measure on $\mathbb{R}$ with no
mass at
$0$ and assume that
%
%
\begin{equation}\label{asso}
n \mathbb{P}(\xi^n_{12} \in\cdot)
\mathop{\longrightarrow}^{w}_{n\to\infty}
\nu\qquad
\mbox{as } n\to\infty.
\end{equation}
Let $G_n$ be the complete network on $\{1,\ldots,n\}$ whose mark on edge
$(i,j)$ equals~$\xi^n_{i j}$. Then, for all integers $B,H$, as $n$
goes to
infinity, in distribution,
\[
(G_n,1) ^{B,H} \longrightarrow(\cT,\eset) ^{B,H}.
\]
Moreover, if $\cT_1,\cT_2$ are independent with common law $\pwit
(\nu)$, then,
in distribution,
\[
( (G_n,1) ^{B,H},(G_n,2) ^{B,H})\longrightarrow( (\cT_1,\eset)
^{B,H}, (\cT_2,\eset) ^{B,H}).
\]
\end{prop}

The second statement is the convergence of the joint law of the finite
networks, where $(G_n,2) ^{B,H}$ is obtained with the same procedure as for
$(G_n,1)^{B,H}$, by starting from the vertex $2$ instead of $1$. In
particular, the second statement implies the first.

This type of convergence is often referred to as \textit{local weak
convergence}, a notion introduced by Benjamini and Schramm
\cite{benjaminischramm} and Aldous and Steele~\cite{aldoussteele}
(see also
Aldous and Lyons~\cite{aldouslyons}). Let us give some examples of application
of this proposition. Consider the case where $\xi^n_{i j} = 1$ with
probability $\lambda/n$ and $\xi^n_{i,j}=\infty$ otherwise. The
network $G_n$
is an Erd\H{o}s--R\'{e}nyi random graph with parameter $\lambda/n$.
From the
proposition, we retrieve the well-known fact that it locally converges
to the
tree of a Yule process of intensity $\lambda$. If $\xi^n_{i,j}= n Y_{i,j}$,
where $Y_{i,j}$ is any nonnegative continuous random variable with density
$1$ at $0+$, then the network converges to $\pwit(\ell_1)$, where
$\ell_1$ is
the Lebesgue measure on $[0,\infty)$. The relevant application for our purpose
is given by the choice
$\xi^n_{i,j}=(a_n / X_{i,j})$, and $\nu=\nu_{\alpha,\theta}$,
where $X_{i,j}$ are
such that $|X_{i,j}|\in\mathbb{H}_\alpha$ and (\ref{theta}) is
satisfied, and
$\nu_{\alpha,\theta}$ is defined by (\ref{appl}). Note that the proposition
applies to all $\alpha>0$ in this setting.
\begin{pf*}{Proof of Proposition~\ref{prop:LWC}}
We order the elements of $J_{B,H}$ in the lexicographic order, that is,
$\eset
\prec1 \prec2 \prec\cdots\prec B \prec11\prec12\prec\cdots\prec B
\cdots
B$. For $\bv\in J_{B,H}$, let $O_{\bv}$ denote the set of offsprings of
$\bv$ in $(G_{n},1)^{B,H}$. By construction, we have $I_\eset= \{1\}$ and
$I_\bv= \sigma_n( \bigcup_{\bw\prec\bv} O_{\bw} )$. At every step
of the
indexing procedure, we sort the marks of the neighboring edges that
have not
been explored at an earlier step $\{1,\ldots, n\} \setminus I_1$,
$\{1,\ldots, n\} \setminus I_{2}, \ldots.$ Therefore, for all
$\mathbf{u}$,
%
%
\begin{equation}\label{eq:lwc1}
\bigl( \xi^n _{ \sigma_n(\mathbf{u}), i} \bigr) _{i \notin I_\mathbf{u}} %
\stackrel{d}{=} %
( \xi^n _{ 1, i} ) _{1 \leq i \leq n - |I _ \mathbf{u}| } .
\end{equation}
Thus, from Lemma~\ref{le:HT} and the independence of the variables
$\xi^n$,
we infer that the marks from a parent to its offsprings in $(G_n,1) ^{B,H}$
converge weakly to those in $(\cT, \eset)^{B,H}$. We now check that all
other marks diverge to infinity. For $\bv, \bw\in J_{B,H}$, we
define
\[
x^n_{\bv,\bw} = \xi^n_{\sigma_n(\bv),\sigma_n(\bw)}.
\]
Also, let
$\{y^n_{\bv,\bw} , \bv, \bw\in J_{B,H}\}$ denote independent variables
distributed as $|\xi^n_{1,2}|$. Let $E^{B,H}$ denote the set of edges
$\{\mathbf{u},\bv\} \in J_{B,H}\times J_{B,H}$ that do not belong to
the finite
tree (i.e., there is no $k\in\{1,\ldots,B\}$ such that $\mathbf
{u}=\bv k$ or
$\bv=\mathbf{u}k$). Lemma~\ref{le:stochdom} below implies that the vector
$\{|x^n_{\bv,\bw}| , \{\bv, \bw\} \in E^{B,H}\}$ stochastically dominates
the vector $\cY^n:=\{y^n_{\bv,\bw} , \{\bv, \bw\} \in E^{B,H}\}
$, that is,
there exists a coupling of the two vectors such that almost surely
$|x^n_{\bv,\bw}| \geq y^n_{\bv,\bw}$, for all $\{\bv, \bw\} \in E^{B,H}$.
Since $J_{B,H}$ is finite (independent of $n$), $\cY^n$ contains a finite
number of variables and (\ref{asso}) implies that the probability of the
event $\{\min_{\{\bv, \bw\} \in E^{B,H}}|x^n_{\bv,\bw}| \leq t\}$
goes to
$0$ as $n\to\infty$, for any \mbox{$t>0$}. Therefore it is now standard to obtain
that if $x_{e}$ denote the mark of edge $e$ in $\cT^{B,H}$, the finite
collection of marks $(x^n_{e})_{e \in J_{B,H} \times J_{B,H} }$
converges in
distribution to $(x_{e})_{e \in J_{B,H} \times J_{B,H}}$ as $n\to
\infty$. In
other words, $(G_n,1) ^{B,H}$ converges in distribution to $(\cT,
\eset)^{B,H}$.

It remains to prove the second statement. It is an extension of the above
argument. We consider the two subnetworks $(G_{n},1)^{B,H}$ and
$(G_{n},2)^{B,H}$ obtained from $(G_n,1)$ and $(G_n,2)$. This gives
rise to
two increasing sequences of sets of vertices $I_{\bv,1}$ and $I_{\bv,2}$
with $\bv\in J_{B,H}$ and two injective maps $\sigma_{n,1}$, $\sigma_{n,2}$
from $J_{B,H}$ to $\{1 , \ldots, n \}$. We need to show that, in
distribution,
%
%
\begin{equation}\label{eq:BH2}
( (G_n,1) ^{B,H},(G_n,2) ^{B,H}) %
\longrightarrow( (\cT_1,\eset) ^{B,H}, (\cT_2,\eset) ^{B,H}).
\end{equation}
Let $V_{n,i}^{B,H} = \sigma_{n,i} ( J_{B,H}) $ be the vertex set of
$(G_{n},i)^{B,H}$, $i=1,2$. There are
\[
C := \frac{B^{H+1} - 1}{B-1}
\]
vertices in $V_{n,i}^{B,H}$, hence the exchangeability of the variables
implies that
\[
\mathbb{P}( 2 \in V_{n,1}^{B,H} ) \leq\frac{C}{n}.
\]
Let $\widetilde G_n = G_n \setminus V_{n,1}^{B,H}$, the subnetwork of $G_n$
spanned by the vertex set $V \setminus V^{B,H}_{n,1}$. Assuming that $2
(B^{H+1} - 1)/(B-1) < n$ and $2 \notin V_{n,1}^{B,H}$, we may then define
$(\widetilde G_{n}, 2)^{B,H}$. If $2 \in V_{n,1}^{B,H}$, $(\widetilde G_{n},
2)^{B,H}$ is defined arbitrarily. The above analysis shows that, in
distribution,
\[
( (G_n,1) ^{B,H},(\widetilde G_n,2) ^{B,H}) %
\longrightarrow( (\cT_1,\eset) ^{B,H}, (\cT_2,\eset) ^{B,H}).
\]
Therefore in order to prove (\ref{eq:BH2}) it is sufficient to prove that
with probability tending to $1$,
\[
V_{n,1}^{B,H} \cap V_{n,2}^{B,H} = \eset.
\]
Indeed, on the event $\{V_{n,1}^{B,H} \cap V_{n,2}^{B,H} = \eset\}$,
$(G_{n},2)^{B,H}$ and $(\widetilde G_{n},2)^{B,H}$ are equal. For $\bv
\in
J_{B,H}$, let $O_{\bv,2}$ denote the set of offsprings of $ \bv$ in
$(G_{n},2)^{B,H}$. We have
\[
I_{\bv,2} = \{2\} \cup\bigcup_{\bw\prec\bv} O_{\bw,2}
\]
and
\begin{eqnarray*}
&&\mathbb{P}( V_{n,1}^{B,H} \cap V_{n,2}^{B,H} \neq\eset) \\
&&\qquad\leq
\mathbb{P}( 2 \in V_{n,1}^{B,H} ) +
\sum_{\bv= \eset}^{B\cdots B}
\mathbb{P}( O_{\bv,2} \cap V_{n,1}^{B,H}
\neq\eset| V_{n,1}^{B,H} \cap I_{\bv,2} = \eset).
\end{eqnarray*}
For any $\mathbf{u}, \bv\in J_{B,H}$, if $V_{n,1}^{B,H} \cap I_{\bv
,2} = \eset$,
then $\sigma_{n,2} ( \bv)$ is neither the ancestor of $\sigma
_{n,1}(\mathbf{u})$,
nor an offspring of $\sigma_{n,1}(\mathbf{u})$. From Lemma \ref
{le:stochdom} below
we deduce that $|\xi^n_{\sigma_{n,1} (\mathbf{u}), \sigma_{n,2}
(\bv)}|$ given
$V_{n,1}^{B,H} \cap I_{\bv,2} = \eset$ dominates stochastically
$|\xi^n_{1,2}|$, and is independent of the i.i.d. vector
$(|\xi^n_{\sigma_{n,2} (\bv),k}|)_{k \in\{1,\ldots,n\}\setminus
(V_{n,1}^{B,H} \cup I_{\bv,2} )}$,
with law $|\xi^n_{1,2}|$. It follows that
\[
\mathbb{P}\bigl(\sigma_{n,1}(\mathbf{u})\in O_{\bv,2} | V_{n,1}^{B,H}
\cap I_{\bv,2} =
\eset\bigr) %
\leq\frac{B}{n - C - | I_{\bv,2} |}.
\]
Therefore,
\begin{eqnarray*}
&&\mathbb{P} ( O_{\bv,2} \cap V_{n,1}^{B,H}
\neq\eset| V_{n,1}^{B,H} \cap I_{\bv,2} = \eset)\\
&&\qquad \leq
\sum_{\mathbf{u}\in J_{B,H} }
\mathbb{P} \bigl( \sigma_{n,1} (\mathbf{u}) \in O_{\bv,2} |
V_{n,1}^{B,H} \cap I_{\bv,2} = \eset\bigr) \\
&&\qquad \leq \frac{ C B }{ n - 2C}.
\end{eqnarray*}
Finally,
\[
\mathbb{P} ( V_{n,1}^{B,H} \cap V_{n,2}^{B,H} \neq\eset) %
\leq\frac C n + \frac{ C^2 B }{ n - 2C},
\]
which converges to $0$ as $n\to\infty$.
\end{pf*}

We have used the following stochastic domination lemma. For any $B,H$
and $n$
let $\cE_n^{H,B}$ denote the (random) set of edges $\{i,j\}$ of the complete
graph on $\{1,\ldots,n\}$, such that
$\{\si_n^{-1}(i),\si_n^{-1}(j)\}$ is not an edge of the finite tree on
$J_{B,H}$. By construction, any loop $\{i,i\}$ belongs to $\cE
_n^{B,H}$. Also,
for $\mathbf{u}\neq\eset$ on the finite tree, let $g(\mathbf{u})$
denote the parent of
$\mathbf{u}$.
\begin{lem}[(Stochastic domination)]\label{le:stochdom}
For any $n\in\mathbb{N}$, and $B,H\in\mathbb{N}$ such that
\[
\frac{B^{H+1} - 1}{B-1} \leq n,
\]
the random variables
\[
\{|\xi^n_{i,j}| , \{i,j\}\in\cE_n^{B,H} \}
\]
stochastically dominate i.i.d. random variables with the same law as law
$|\xi^n_{1,2}|$. Moreover, for every $\eset\neq\mathbf{u}\in
J_{B,H}$, the random
variables
\[
\bigl\{\bigl|\xi^n_{\si_n(\mathbf{u}),i}\bigr| , i\in\{1,\ldots,n\}\setminus\si
_n(g(\mathbf{u}))\bigr\},
\]
stochastically dominate i.i.d. random variables with the same law as law~$|\xi^n_{1,2}|$.
\end{lem}

\begin{pf}
The censoring process which deletes the edges that belong to the tree on
$J_{B,H}$ has the property that at each step the $B$ lowest absolute values
are deleted from some \textit{fresh} (previously unexplored) subset of edge
marks. Using this and the fact that the edge marks $\xi^n_{i,j}$ are
i.i.d. we see that both claims in the lemma are implied by the
following simple
statement.

Let $Y_1,\ldots,Y_m$ denote i.i.d. positive random variables. Suppose $m=n_1
+ \cdots+ n_\ell$, for some positive integers $\ell$, $n_1,\ldots
,n_\ell$, and
partition the $m$ variables in $\ell$ blocks $I^1,\ldots,I^\ell$ of
$n_1,\ldots,n_\ell$ variables each. Fix some nonnegative integers $k_j$
such that $k_j\leq n_j$ and call $q_1^j,\ldots,q^j_{k_j}$, the (random)
indexes of the $k_j$ lowest values of the variables in the block $I^j$ (so
that $Y_{q^1_1}$ is the lowest of the $Y_1,\ldots,Y_{n_1}$, $Y_{q^1_2}$ is
the second lowest of the $Y_1,\ldots,Y_{n_1}$ and so on). Consider the random
index sets of the $k_j$ minimal values in the $j$th block,
$J^j:=\bigcup_{i=1}^{k_j}\{q^j_i\}$, and set $J=\bigcup_{j=1}^\ell
J^j$. If
$k_j=0$ we set $J^j=\eset$. Finally, let $\wt Y$ denote the vector $\{
Y_i ,
i=1,\ldots,m ; i\notin J\}$. Then we claim that $\wt Y$ stochastically
dominates $m-\sum_{j=1}^\ell k_j$ i.i.d. copies of $Y_1$.\vspace*{1pt}

Indeed, the coupling can be constructed as follows. We first extract a
realization $y_1,\ldots,y_m$ of the whole vector. Given this we isolate the
index sets $J^1,\ldots,J^\ell$ within each block. We then consider two
vectors $\cZ,\cV$ obtained as follows. The vector
$\cZ_1=(z^1_1,\ldots,z^1_{n_1-k_1},z^2_1,\ldots,z^2_{n_2-k_2},\ldots
,z^\ell_{n_\ell-k_\ell})$
is obtained by extracting the $n_1-k_1$ values $z^1_1,\ldots,z^1_{n_1-k_1}$
uniformly at random (without replacement) from the values
$y_1,\ldots,y_{n_1}$ (in the block $I^1$), the $n_2-k_2$ variables
$z^2_1,\ldots,z^2_{n_2-k_2}$ in the same way from the values
$y_{n_1+1},\ldots,y_{n_1+n_2}$ (in the block\vspace*{1pt} $I^2$), and so on. On the other
hand, the vector $\cV%
=(v^1_1,\ldots,v^1_{n_1-k_1},v^2_1,\ldots,v^2_{n_2-k_2},\ldots,v^\ell
_{n_\ell-k_\ell})$
is obtained as follows. For the first block we take $v^1_i$, $i=1,\ldots,
n_1-k_1$ equal to $z^1_i$ whenever an index $i\in I^1\setminus J^1$ was
picked for the vector $z^1_1,\ldots,z^1_{n_1-k_1}$, and we assign the
remaining values (if any) through an independent uniform permutation of
those variables $y_i, i\in I^1\setminus J^1$ which were not picked for the
vector $z^1_1,\ldots,z^1_{n_1-k_1}$. We repeat this procedure for all other
blocks to assign all values of $\cV$. By construction, $\cV\geq\cZ$
coordinate-wise. The conclusion follows from the observation that $\cZ
$ is
distributed like a vector of $m-\sum_{j=1}^\ell k_j$ i.i.d. copies of
$Y_1$, while $\cV$ is distributed like our vector $\wt Y$.
\end{pf}

\subsection{\texorpdfstring{Proof of Theorem \protect\ref{th:convop}}{Proof of Theorem 2.3}}

\mbox{}

\begin{pf*}{Proof of Theorem~\ref{th:convop}\textup{(i)}}
Let $\nu=\nu_{\alpha,\theta}$ be as in (\ref{appl}), and let
$(\cT_\alpha,\eset)$ be a realization of the $\pwit(\nu)$. The
mark on edge
$(\bv,\bv k)$ in $\cT_\alpha$ is denoted by $x_{(\bv, \bv k)}$ or simply
$x_{\bv k}$. By definition, we have $x_{(\bv,\bw)} = \infty$ if $\bv
$ and
$\bw$ are at graph-distance different from $1$. In particular, if we set
$y_\bv= \sign(x_{\bv})| x_{\bv}|^{\alpha}$, then the point sets
$\Xi_\bv=
\{y_{\bv k } \}_{k \geq1}$ are independent Poisson point processes of
intensity\vadjust{\goodbreak} $\ell_\theta= \theta\ind_{\{x > 0\}}\,dx + (1-\theta) \ind
_{\{x <
0\}}\,dx$. We may thus build a realization of the operator $\bT$ on
$\cT_\alpha$ [cf. (\ref{tone})]. Let $G_n$ be the complete network on
$\{1,\ldots,n\}$ whose mark on edge $(i,j)$ is $\xi^n_{i,j}:=a_n/X_{i,j}$.
Next,\vspace*{-1pt} we apply Proposition~\ref{prop:LWC}. For all $B$, $H$, $(G_n,1)^{B,H}$
converges weakly to $(\cT_\alpha,\eset)^{B,H}$. Let $\si_n^{B,H}$
be the map
$\si_n$ associated with the network $(G_n,1)^{B,H}$ (see the construction
given before Proposition~\ref{prop:LWC}). From the Skorokhod representation
theorem we may assume that $(G_n,1)^{B,H}$ converges a.s. to
$(\cT_\alpha,\eset)^{B,H}$ for all $B,H$. Thus we may find sequences
$B_n,H_n$ tending to infinity, such that $(B_n ^{H_n +1} - 1 )/ (B_n - 1)
\leq n$ and such that for any pair $\mathbf{u},\bv\in\mathbb{N}^f$
we have
$\xi^n_{(\wt\sigma_n (\mathbf{u}), \wt\sigma_n (\bv))} \to
x_{(\mathbf{u}, \bv)}$ a.s. as
$n\to\infty$, where $\wt\si_n:=\si_n^{B_n,H_n}$. The map $\wt\si
_n$ can be
extended to a bijection $\mathbb{N}^f\to\mathbb{N}^f$. It follows
that a.s. %
%
%
\begin{equation}\label{pointw}
\langle
\delta_\mathbf{u}, \wt\sigma_n^{-1} (a_n^{-1} X ) \wt\sigma_n
\delta_\bv\rangle
= \frac{ 1}{ \xi^n_{(\wt\sigma_n (\mathbf{u}),\wt\sigma_n (\bv
))} } \to\frac{
1}{x_{(\mathbf{u}, \bv)} }= \langle\delta_\mathbf{u}, \bT\delta
_\bv\rangle.
\end{equation}
Fix\vspace*{1pt} $\bv\in\mathbb{N}^f$, and set $\psi_n^\bv:= \wt\si_n^{-1} (a_n^{-1}
X)\wt\si_n \delta_\bv$. To prove Theorem~\ref{th:convop}(i) it is sufficient
to show
that 
$\psi_n^\bv\to\bT\delta_\bv$ in $L^2( \mathbb{N}^f)$ almost
surely as
$n\to\infty$, that is,
%
%
\begin{equation}\label{l2con}
\sum_{\mathbf{u}}
(\langle\delta_\mathbf{u},\psi_n^\bv\rangle-
\langle\delta_\mathbf{u},\bT\delta_\bv\rangle)^2
\to0 .
\end{equation}
Since from (\ref{pointw}) we know that $\langle
\delta_\mathbf{u},\psi_n^\bv\rangle\to\langle
\delta_\mathbf{u},\bT\delta_\bv\rangle$ for every $\mathbf
{u}$, the claim follows
if we have (almost surely) uniform (in $n$) square-integrability of
$(\langle\delta_\mathbf{u},\psi_n^\bv\rangle
)_\mathbf{u}$. This in turn
follows from Lemmas~\ref{le:stochdom} and~\ref{le:PoiExt}(i).
The proof
of Theorem~\ref{th:convop}(i) is complete.
\end{pf*}
\begin{pf*}{Proof of Theorem~\ref{th:convop}\textup{(ii)}}
We need the following two facts:
%
%
\begin{equation}\label{eq:llnweak}
\lim_{n\to\infty}\frac{\rho_1}{n}=1 \qquad\mbox{in probability} ,
\end{equation}
and there exists $\delta>0$ such that
%
%
\begin{equation}\label{eq:llninf}
\liminf_{n\to\infty}\min_{1\leq i\leq n}\frac{\rho_i}{n}>\delta\qquad
\mbox{a.s.}
\end{equation}
Clearly, (\ref{eq:llnweak}) is a law of large numbers and holds actually
a.s. (recall that for $\alpha>1$ we assume the mean of $U_{i,j}$ to
be $1$).
Let us establish the a.s. uniform bound (\ref{eq:llninf}). For every
$\epsilon>0$, there exists $R>0$ such that
$\mathbb{E}(U_{i,j}\ind_{\{U_{i,j}<R\}})\geq1-\epsilon$. If we define
$\rho_i^R=\sum_{j=1}^nU_{i,j}\ind_{\{U_{i,j}<R\}}$, then
%
\[
\liminf_{n\to\infty}\min_{1\leq i\leq n}\frac{\rho_i}{n}
\geq\liminf_{n\to\infty}\min_{1\leq i\leq n}\frac{\rho_i^R}{n} .
\]
Therefore (\ref{eq:llninf}) is implied by the uniform law of large numbers
in~\cite{bordenave-caputo-chafai}, Lem\-ma~2.2, 
applied to the bounded variables $U_{i,j}\ind_{\{U_{i,j}<R\}}$.

Next, we claim that for all $\mathbf{u}\in\mathbb{N}^f$, in probability
%
%
\begin{equation}\label{eq:sigmalln}
\lim_{n\to\infty} \frac{\rho_{\wt\sigma_n( \mathbf{u}) }}{n} =1.\vadjust{\goodbreak}
\end{equation}
To prove this we first observe that by Lemma~\ref{le:stochdom} and
(\ref{eq:llnweak}) we have in probability
\[
\limsup_{n\to\infty} \frac{ ( \rho_{\wt\sigma_n( \mathbf
{u}) } -
U_{\wt\sigma_n( \mathbf{u}) , \wt\sigma_n(g (\mathbf{u}))}
) } {n} \leq1.
\]
On the other hand $U_{\wt\sigma_n( \mathbf{u}) , \wt\sigma_n(g
(\mathbf{u}))}$ is
stochastically dominated by the maximum of $n$ i.i.d. variables with law
$U_{i,j}$. The latter converges in distribution on the scale $a_n$
[cf. Lemma
\ref{le:PoiExt}(i)], and we know that $a_n/n \to0$. It follows that in
probability $\limsup_{n\to\infty} \rho_{\wt\sigma_n( \mathbf{u})
}/ n \leq1$.
Next, we can estimate
\[
\rho_{\wt\sigma_n( \mathbf{u}) } \geq
\sum_{i\in\{1,\ldots, n\}\setminus I_\mathbf{u}} U_{\wt\sigma_n(
\mathbf{u}),i} .
\]
Now, observe that if $\mathbf{u}\in\mathbb{N}^f$ belongs to
generation $h$, then
the set $I_\mathbf{u}$ contains at most $O(B_n^h)$ elements, while $n$
is at least
of order $B_n^{H_n}$, where $B_n,H_n$ are the sequences used in the
proof of
Theorem~\ref{th:convop}(i). In particular, it follows that
$|I_\mathbf{u}|=o(n)$
and therefore (\ref{eq:lwc1}) and (\ref{eq:llnweak}) imply that
$\liminf_{n\to\infty} \rho_{\wt\sigma_n( \mathbf{u}) }/ n\geq1$
in probability.
This proves (\ref{eq:sigmalln}).

Thanks to (\ref{eq:sigmalln}), from the Slutsky lemma and the Skorokhod
representation theorem, we may also assume that for each
$\bv\in\mathbb{N}^f$, $\rho_{\wt\sigma_n(\bv)}/n$ converges
a.s. to $1$. We
need to show that for each $\bv\in\mathbb{N}^f$, (\ref{l2con})
holds with
the new vector $\psi_n^\bv:= \wt\si_n^{-1} (\kappa_n S)\wt\si_n
\delta_\bv$, 
\[
\langle\delta_\bw,\psi_n^\bv\rangle= \kappa_n
\frac{U_{\wt\si_n (\bw),\wt\si_n(\bv)}}{\sqrt{\rho_{\wt\si_n
(\bv)}\rho_{\wt\si_n (\bw)}}} .
\]
Thanks to (\ref{eq:llninf}), $(\langle
\delta_\bw,\psi_n^\bv\rangle)_\bw$ is uniformly square-integrable
[cf. the proof of (\ref{l2con})], and all we have to check is that
$(\langle\delta_\bw,\psi_n^\bv\rangle- \langle
\delta_\bw,\bT\delta_\bv\rangle)^2\to0$ for fixed $\bw$.
Here $\bT$
is the operator appearing in the proof of Theorem~\ref{th:convop}(i) above,
now with the choice $\theta=1$. We have
\begin{eqnarray*}
&&(\langle
\delta_\bw,\psi_n^\bv\rangle- \langle
\delta_\bw,\bT\delta_\bv\rangle)^2
\\
&&\qquad \leq2 \bigl(a_n^{-1}U_{\wt\si_n (\bw),\wt\si_n (\bv
)}\bigl(1-n/
\sqrt{\rho_{\wt\si_n (\bv)}\rho_{\wt\si_n (\bw)}}
\bigr)\bigr)^2\\
&&\qquad\quad{}+ 2\bigl(a_n^{-1}U_{\wt\si_n (\bw),\wt\si_n (\bv)} - \langle
\delta_\bw,\bT\delta_\bv\rangle\bigr)^2 .
\end{eqnarray*}
The second term above converges to zero as in the proof of point (i). For
the first term we use $\rho_{\wt\sigma_n(\bv)}/n \to1$ and
$\rho_{\wt\sigma_n(\bw)}/n \to1$. This proves point (ii).
\end{pf*}
\begin{pf*}{Proof of Theorem~\ref{th:convop}\textup{(iii)}}
The setting is as in the proof of point~(ii) above, but now
$\alpha\in(0,1)$. We build the operator $\bS$ on the tree $\cT
_\alpha$ as in
(\ref{opsym}). We need to prove that for any $\bv\in\mathbb{N}^f$,
a.s.
%
%
\begin{equation}
\label{tos}
\sum_{\bw}(\langle
\delta_\bw,\psi_n^\bv\rangle-
\langle
\delta_\bw,\bS\delta_\bv\rangle)^2 \to0 ,
\end{equation}
with $\psi_n^\bv:= \wt\si_n^{-1} S \wt\si_n \delta_\bv$, that is,
\[
\langle\delta_\bw,\psi_n^\bv\rangle= \frac{U_{\wt
\si_n
(\bw),\wt\si_n(\bv)}}{\sqrt{\rho_{\wt\si_n (\bv)}\rho_{\wt
\si_n
(\bw)}}} .
\]
Let us first show that for any $\bv,\bw\in\mathbb{N}^f$ we have
a.s.
%
%
\begin{equation}
\label{eq:convK}
\frac{U_{\wt\sigma_n (\bw), \wt\sigma_n (\bv)} }
{
\sqrt{\rho_{\wt\si_n (\bv)}\rho_{\wt\si_n (\bw)}}
}
\to\frac{\langle
\delta_\bw,
\bT\delta_\bv\rangle}
{\sqrt{\rho(\bv)\rho(\bw)}} = \langle
\delta_\bw,\bS\delta_\bv\rangle.
\end{equation}
Multiplying and dividing by $a_n$ and using (\ref{pointw}) with
$\theta=1$,
we see that (\ref{eq:convK}) holds if
%
%
\begin{equation}\label{rhoss}
a_n^{-1}\rho_{\wt\si_n (\bv)} \to\rho(\bv) ,
\end{equation}
almost surely, for every $\bv\in\mathbb{N}^f$. In turn, (\ref
{rhoss}) can be
proved as follows. Let \mbox{$k \in\mathbb{N}$}, and consider the tree with vertex
set $J_{k,k}$, obtained as in Proposition~\ref{prop:LWC} with $B=H=k$. Since
$J_{k,k}$ is a finite set, for any $\bv$, (\ref{pointw}) implies that a.s.
\[
a_n^{-1}\sum_{\mathbf{u}\in J_{k,k}} U_{\wt\sigma_n (\bv), \wt
\sigma_n (\mathbf{u}) }\to
{\sum_{\mathbf{u}\in J_{k,k}} x^{-1}_{\bv,\mathbf{u}}} .
\]
By Lemmas~\ref{le:stochdom} and~\ref{le:PoiExt}(ii), $\sum
_{\mathbf{u}\notin
J_{k,k}} a_n^{-1} U_{\wt\sigma_n (\bv), \wt\sigma_n (\mathbf
{u})}$ a.s. converges uniformly (in~$n$) to $0$ as $k$ goes to
infinity. This proves
(\ref{eq:convK}) and (\ref{rhoss}).

Once we have (\ref{eq:convK}), to conclude the proof it is sufficient to
show that a.s.
%
%
\begin{equation}\label{unifint} \lim_{k\to\infty} \sup_n
\sum_{\bw\notin J_{k,k}}(\langle\delta_\bw,\psi_n^\bv
\rangle)^2
= 0 .
\end{equation}
However,\vspace*{1pt} using (\ref{rhoss}) and the simple bound $(\langle
\delta_\bw,\psi_n^\bv\rangle)^2\leq\frac{U_{\wt\si_n
(\bv),\wt\si_n(\bw)}}{\rho_{\wt\si_n (\bv)}}$, we have that
(\ref{unifint}) again follows from an application of Lemmas~\ref{le:stochdom}
and~\ref{le:PoiExt}(ii). This completes the proof of Theorem
\ref{th:convop}(iii).
\end{pf*}

\subsection{Two-points local operator convergence}

In the proof of the main theorems, we will need a stronger version of Theorem
\ref{th:convop}. Define the $2 n \times2 n$ matrices
\[
X \oplus X \quad\mbox{and}\quad S \oplus S,
\]
where ``$\oplus$'' denotes the usual direct sum decomposition,
$X\oplus X
(\phi_1,\phi_2) = (X \phi_1$, $X \phi_2)$, for $n$-dimensional vectors
$\phi_1,\phi_2$. As for the limiting operators, we realize them on
the Hilbert
space $L^2(V)\oplus L^2(V)$ with $V=\mathbb{N}^f$. We consider two independent
realizations $\cT^1_\alpha$, $\cT^2_\alpha$ of the PWIT($\ell
_\theta$), and
call $\bT_1, \bS_1,\bT_2,\bS_2$ the associated operators as in Section
\ref{limop}. We may then define
\[
\bT_1 \oplus\bT_2 \quad\mbox{and}\quad \bS_1\oplus\bS_2.
\]
By Proposition~\ref{prop:LWC}, $((G_n,1))^{B,H},(G_n,2)^{B,H})$ converges
weakly to $((\cT^1_\alpha, \eset)^{B,H} $, $(\cT^2_\alpha, \eset
)^{B,H} )$. As
before we can view the matrices $X\oplus X$ and $S\oplus S$ as bounded
self-adjoint operators on $L^2(V)\oplus L^2(V)$. Therefore, arguing as
in the
proof of Theorem~\ref{th:convop}, it follows that, in distribution,
for all
$(\phi_1, \phi_2) \in\cD\times\cD$,
\[
\sigma_{n}^{-1} a_n^{-1} X \oplus X \sigma_n (\phi_1,\phi_2)
\rightarrow\bT_1\oplus\bT_2 (\phi_1,\phi_2) ,
\]
where, $\sigma_n = \sigma_n^1 \oplus\sigma_n^2$, and, as above, for
$i \in
\{1,2\}$, $\sigma_n^i$ is a bijection on $\mathbb{N}^f$, extension of the
injective indexing map from $\mathbb{N}^f$ to $\{1,\ldots,n\}$, such that
$\sigma_n^i (\eset) = i$. Analogous convergence results hold for the matrix
$S\oplus S$. We can thus extend the statement of Theorem \ref
{th:convop} to
the following local convergence of operators in $L^2(V)\oplus L^2(V)$. To
avoid lengthy repetitions we omit the details of the proof.
\begin{theorem}\label{th:convop2}
As $n$ goes to infinity, in distribution:
\begin{enumerate}[(iii)]
\item[(i)] if $\alpha\in(0,2)$, then %
$(a_n ^{-1} X\oplus a_n^{-1} X ,(1,2) ) \to(\bT_1 \oplus\bT_2
,(\eset,\eset))$;
\item[(ii)] if $\alpha\in(1,2)$ and $\theta=1$, then %
$( \kappa_n S \oplus\kappa_n S ,(1,2)) \to(\bT_1 \oplus\bT_2 ,
(\eset,\eset) )$;
\item[(iii)] if $\alpha\in(0,1)$, then %
$(S\oplus S ,(1,2)) \to(\bS_1 \oplus\bS_2 ,(\eset,\eset))$.
\end{enumerate}
\end{theorem}

\section{Convergence of the empirical spectral distributions}
\label{se:convesd}

\subsection{\texorpdfstring{Markov matrix, $\alpha\in(0,1)$: Proof of Theorem \protect\ref{th:k01}}
{Markov matrix, alpha in (0,1): Proof of Theorem 1.4}}
\label{sec:k01}

Recall that $\bS$ is a bounded self-adjoint operator on $L^2(V)$, whose
spectrum is contained in $[-1,1]$ [cf. (\ref{isom})]. The resolvents
of $S$
and $\bS$ are the functions on $\mathbb{C}_+ = \{z \in\mathbb{C} \dvtx
\Im z > 0
\}$:
\[
R^{(n)} (z) %
= ( S - z I) ^{-1} \quad\mbox{and}\quad R (z) = ( \bS- z I) ^{-1}.
\]
For $\ell\in\mathbb{N}$, set
%
%
\begin{equation}\label{gammal}
\mathbf{p}_\ell:= \langle\delta_\eset, \bS^{\ell}\delta_\eset
\rangle.
\end{equation}
Note that $\mathbf{p}_\ell=\frac1{\rho(\eset)} \langle\delta_\eset
, \bK^{\ell}
\delta_\eset\rangle_\rho$ is the probability that the random walk
on the PWIT
associated with the stochastic operator $\bK$ comes back to the root
(where it
started) after $\ell$ steps. In particular, $\mathbf{p}_\ell=0$ for
$\ell$ odd.
We set $\mathbf{p}_0=1$. Let $\mu_\eset$ denote the spectral measure of
$\bS$
associated with $\delta_\eset$ (see e.g.,~\cite{reedsimon}, Chapter VII).
Equivalently, $\mu_\eset$ is the spectral measure of $\bK$
associated with the
$L^2(V,\rho)$ normalized vector $\hat\delta_\eset:=
\delta_\eset/\sqrt{\rho(\eset)}$ [cf. (\ref{isom})]. In particular,
$\mu_\eset$ is a probability measure supported on $[-1,1]$ and such
that $\mathbf{p}_\ell= \int_{-1}^{1}x^\ell\mu_\eset(dx)$, for every
$\ell$.
Since all odd
moments vanish $\mu_\eset$ is symmetric. Moreover, for any $z\in
\mathbb{C}_+$
we have
\[
\langle\delta_\eset, R (z)
\delta_\eset\rangle= \int_{-1}^1 \frac{ \mu_\eset(dx)}{x - z},
\]
that is, $\langle\delta_\eset, R (z)\delta_\eset\rangle$ is the
Cauchy--Stieltjes transform of $\mu_\eset$. Recall that the Cauchy--Stieltjes
transform\vadjust{\goodbreak} of a probability measure $\mu$ on $\mathbb{R}$ is the analytic
function on $\mathbb{C}_+$ given by
\[
m_\mu(z) = \int_\mathbb{R} \frac{\mu(dx)}{x - z}.
\]
The function $m_\mu$ characterizes the measure $\mu$, $|m_\mu(z) |
\leq(\Im
z) ^{-1}$, and weak convergence of $\mu_n$ to $\mu$ is equivalent to the
convergence $m_{\mu_n}(z)\to m_\mu(z)$ for all $z\in\mathbb{C}_+$. By
construction
\[
\frac1 n \tr R^{(n)} (z) %
= \int_{-1}^1 \frac{ \mu_K (dx)}{x - z} %
= m_{\mu_K} (z) ,
\]
where $\mu_K$ is the ESD of $K$, which coincides with the ESD of $S$.
Using exchangeability and linearity, we get
\[
\mathbb{E} R^{(n)}_{1,1} (z)
= \mathbb{E} m_{\mu_K} (z) = m_{\mathbb{E} \mu_K} (z).
\]
Since $R^{(n)}(z) _{1,1} \leq(\Im z) ^{-1}$ is bounded, we may apply
Theorems
\ref{th:strongres} and~\ref{th:convop}, and obtain, for all
$z \in
\mathbb{C}_+$,
%
%
\begin{equation}\label{eq:Kav}
\lim_{n \to\infty} m_{\mathbb{E} \mu_K} (z) = m_{\mathbb{E} \mu
_\eset} (z).
\end{equation}
We define
\[
\wt\mu_\alpha= \mathbb{E} \mu_\eset.
\]
%
Next, we shall prove that, for all $z \in\mathbb{C}_+$,
%
%
\begin{equation}\label{eq:KL1}
\lim_{n \to\infty} %
\mathbb{E} | m_{ \mu_K} (z) - m_{\mathbb{E} \mu_\eset} (z)
| = 0.
\end{equation}
We have
\[
\mathbb{E} | m_{\mu_K} (z) - m_{ \mathbb{E} \mu_\eset} (z) | %
\leq\mathbb{E} | m_{\mu_K} (z) - \mathbb{E} m_{\mu_K} (z)
| %
+ | m_{\mathbb{E} \mu_K} (z) - m_{\mathbb{E} \mu_\eset}
(z)|.
\]
On the right-hand side, the second term converges to $0$ by (\ref{eq:Kav}).
The first term is equal to
\[
\mathbb{E} \Biggl| \frac1n\sum_{k=1}^n \bigl[ R ^{(n)}_{k,k} (z)-
\mathbb{E}
R ^{(n)}_{k,k} (z) \bigr] \Biggr| .
\]
By exchangeability, we note that
\begin{eqnarray*}
&&\mathbb{E} \Biggl[\Biggl( \frac1 n \sum_{k=1}^n \bigl[ R
^{(n)}_{k,k} (z)-
\mathbb{E}
R ^{(n)}_{k,k} (z) \bigr] \Biggr)^2\Biggr] \\
&&\qquad = \frac1 n \mathbb{E} \bigl(R ^{(n)}_{1,1} - \mathbb{E} R
^{(n)}_{1,1} \bigr)^2 + \frac{n(n-1)}{n^2} \mathbb{E} \bigl[ \bigl(R
^{(n)}_{1,1} - \mathbb{E} R ^{(n)}_{1,1} \bigr)\bigl(R
^{(n)}_{2,2} - \mathbb{E} R ^{(n)}_{2,2}\bigr) \bigr] \\
&&\qquad \leq\frac{1}{n (\Im z)^2} + \mathbb{E} \bigl[ \bigl(R
^{(n) }_{1,1}
- \mathbb{E} R ^{(n)}_{1,1} \bigr)\bigl(R ^{(n) }_{2,2} - \mathbb
{E} R
^{(n) }_{2,2}\bigr) \bigr].
\end{eqnarray*}
Theorems~\ref{th:strongres} and~\ref{th:convop2} imply that
$(R_{1,1}(z),R_{2,2}(z))$ are asymptotically independent. Since these
variables are bounded, they are also asymptotically uncorrelated, and
(\ref{eq:KL1}) follows.\vadjust{\goodbreak}

Finally, observe that the sequence of measures $\mu_K$ is a.s. tight.
Therefore the convergence (\ref{eq:KL1}) is sufficient to establish
a.s. convergence of $\mu_K$ to $\wt\mu_\alpha$. This completes the proof
of Theorem
\ref{th:k01}.

\subsection{\texorpdfstring{I.i.d. matrix, $\alpha\in(0,2)$: Proof of Theorem \protect\ref{th:iida}}
{I.i.d. matrix, alpha in (0,2): Proof of Theorem 1.2}}\label{sec:iida}

Set $A_n = a_n^{-1} X$. For $z\in\mathbb{C}_+$, we define the
Cauchy--Stieltjes transform,
\[
m_{A_n} (z) = \int\frac{d \mu_{A_n}(x)} {x - z}
= \frac1 n \sum_{k=1}^n R^{(n)}_{k,k} (z) ,
\]
where
\[
R^{(n)} (z) = (A_n - z I)^{-1} ,
\]
is the resolvent of $A_n$. By exchangeability, $\mathbb{E} m_{A_n} (z) =
\mathbb{E} R^{(n)}_{1,1}(z)$. From Proposition~\ref{esssa} we know
that $\bT$
is self-adjoint. Therefore from Theorems~\ref{th:strongres} and~\ref{th:convop} we infer
%
%
\begin{equation}\label{meancon}
\mathbb{E} m_{A_n} (z) %
\to\mathbb{E} h(z) ,\qquad %
h(z):=\langle\delta_\eset, (\bT-z I)^{-1}\delta_\eset\rangle.
\end{equation}
As in the proof of Theorem~\ref{th:k01} we may write $\mathbb{E} h(z) =
\mathbb{E} m_{\mu_\eset}=m_{\mathbb{E}\mu_\eset}$, that is the
Cauchy--Stieltjes transform of the expected value of the random spectral
measure $\mu_\eset$ associated to $\bT$ at the root vector $\delta
_\eset$.
From (\ref{meancon}) we obtain the weak convergence of $\mathbb{E}
\mu_{A_n}$
to $\mu_\alpha:=\mathbb{E}\mu_\eset$. To obtain a.s. weak
convergence of
$\mu_{A_n}$ to $\mu_\alpha$, from Lemma~\ref{astight} it suffices
to prove the
$L^1$ convergence of Cauchy--Stieltjes transforms as in (\ref
{eq:KL1}). This
in turn is obtained by repeating word by word the argument in the proof of
Theorem~\ref{th:k01}.

Thus, we have obtained $\mu_{A_n}\to\mu_\alpha$ almost surely.
Since the
operator $\bT$ only depends on the two parameters $\alpha$ and
$\theta$, where
the latter is defined by (\ref{theta}), the LSD $\mu_\alpha$ might still
depend on the parameter $\theta$. However, the fact that $\mu_\alpha
$ is
independent of $\theta$ follows from Lemma~\ref{le:itunique} below, which
implies in particular that the values $m_{\mu_\alpha}(i t)=\mathbb{E}[h(i
t)]$, $t>0$, are uniquely determined by $\alpha$, and therefore by
analyticity, all values $m_{\mu_\alpha}(z)$, $z\in\mathbb{C}_+$ are uniquely
determined by $\alpha$. This ends the proof of Theorem~\ref{th:iida}.

We remark that in the proof of Theorem~\ref{th:iida} one can avoid
establishing (\ref{eq:KL1}) plus almost sure tightness [Lemma
\ref{astight}(i)] as we do above. Namely, the convergence of
expected values
$\mathbb{E}\mu_{a_n^{-1}X}\to\mu_\alpha$ is sufficient. This\vspace*{1pt}
follows from an
a priori concentration estimate (see \cite
{bordenave-caputo-chafai-heavygirko}).
However, we did that piece of extra work here since we need it anyway
in the
case of Markov matrices, where the mentioned concentration estimate is not
available.


\subsection{\texorpdfstring{Markov matrix, $\alpha\in(1,2)$: Proof of Theorem \protect\ref{th:k12}}
{Markov matrix, alpha in (1,2): Proof of Theorem 1.3}}
\label{sec:k12}

The proof given above for the matrix $A_n=a_n^{-1}X$ applies without
modifications to the new matrix $A_n:=\kappa_n S$, where $ S_{i,j} =
\frac{U_{i,j}}{\sqrt{\rho_i \rho_j}}$. In particular, we use Theorems
\ref{th:convop}(ii), \ref{th:convop2}(ii) and
Lem\-ma~\ref{astight}(ii) to obtain the a.s. weak convergence of $\mu
_{A_n}$ to
$\mu_\alpha=\mathbb{E}\mu_\eset$, where $ \mu_\eset$ is the
random spectral
measure of $\bT$ at the root. This ends the proof of Theorem~\ref{th:k12}.

\subsection{\texorpdfstring{Markov matrix, $\alpha=1$: Proof of Theorem \protect\ref{th:k1}}
{Markov matrix, alpha = 1: Proof of Theorem 1.5}}
\label{sec:k1}

Suppose now that $\alpha=1$ and set $w_n=\int_0^{a_n} x\cL(dx)$ and
$\kappa_n=na_n^{-1}w_n$. A close inspection of the proof of Theorem
\ref{th:convop}(ii) and Theorem~\ref{th:k12} reveals that all
arguments used
for $\alpha\in(1,2)$ can be applied to the case $\alpha=1$ without
modifications except for the two estimates (\ref{eq:llnweak}) and
(\ref{eq:llninf}), which have to be replaced by (\ref{eq:llnweaka1}) and
(\ref{eq:llninfa1}) below, respectively. For (\ref{eq:llninfa1}) we
shall use
the hypothesis (\ref{eq:conda1wn}) on~$w_n$. Let us start by proving
that, in
probability
%
%
\begin{equation}\label{eq:llnweaka1}
\lim_{n\to\infty} \frac{\rho_1}{n w_n } = 1.
\end{equation}
We recall that, for fixed $i$, $a_n^{-1}(\rho_i - n w_n)$ converges in
distribution to a $1$-stable law (see, e.g.,~\cite{zinn81},
Theorem 1).
Therefore it suffices to show that $\kappa_n=a_n^{-1} n w_n \to\infty
$. To see
this we may argue as follows. Observe that, for any $\varepsilon>0$
\[
\kappa_n = \mathbb{E}\sum_{i=1}^n a_n^{-1} V_i \ind_{\{
a_n^{-1}V_i\leq
1\}} %
\geq\mathbb{E}\sum_{i=1}^n a_n^{-1} V_i \ind_{\{\varepsilon\leq
a_n^{-1}V_i\leq1\}},
\]
where $V_1\geq V_2\geq\cdots$ are the ranked values of $U_{1,j}$,
$j=1,\ldots,
n$. From Lem\-ma~\ref{le:PoiExt}(i) the right-hand\vspace*{1pt} side above, for any
$\varepsilon>0$, converges to $\mathbb{E}\sum_{i} x_i \ind_{\{
\varepsilon\leq
x_i\leq1\}}$, where the $x_i$ are distributed according to the PPP with
intensity $x^{-2}\,dx$ on $(0,\infty)$. While this sum is finite for every
$\varepsilon>0$ it is easily seen to diverge (logarithmically) for
$\varepsilon\to0$. This achieves the proof of (\ref{eq:llnweaka1}).

Next, we claim that if $w_n$ satisfies (\ref{eq:conda1wn}), then there exists
$\delta>0$ such that, a.s.
%
%
\begin{equation}\label{eq:llninfa1}
\liminf_{n\to\infty} \min_{1 \leq i \leq n} \frac{\rho_i}{n w_n}
\geq\delta.
\end{equation}
To establish (\ref{eq:llninfa1}), let us define $b_n = a_{\lfloor
n^\varepsilon\rfloor}$ so that $\mathbb{E}(U_{1,i} \ind_{\{U_{1,i}
\leq b_n\}}) =
w_{\lfloor n^\varepsilon\rfloor} $ and
\[
\rho_1 \geq S_n := \sum_{i = 1}^n U_{1,i}\ind_{\{U_{1,i}\leq b_n\}}.
\]
From the union bound,
\[
\mathbb{P}\biggl(\min_{1 \leq i \leq n}\frac{\rho_i}{n w_n}<\delta
\biggr) %
\leq n \mathbb{P}\biggl(\frac{\rho_1}{n w_n} < \delta\biggr).
\]
From the Borel--Cantelli lemma, it is thus sufficient to prove that for some
$\delta>0$
%
%
\begin{equation}\label{eq:bcSn}
\sum_{n \geq1} n \mathbb{P}(S_n < \delta n w_n) <
\infty.
\end{equation}
By assumption, there exists $\delta> 0$ such that for all $n$ large enough,
$w_{\lfloor n^\varepsilon\rfloor} \geq2 \delta w_n$. We define
\[
V_i = U_{i,1} \ind_{\{U_{1,i} \leq b_n\}} - w_{\lfloor n^\varepsilon
\rfloor} %
\quad\mbox{and}\quad %
\overline S_n = \sum_{i = 1} ^ n V_i.\vadjust{\goodbreak}
\]
Note that $\mathbb{E} V_i = \mathbb{E} \overline S_n = 0$. We get for
all $n$
large enough
%
%
\begin{equation} \label{eq:Sn} %
\mathbb{P}(S_n < \delta n w_n) %
= \mathbb{P}\bigl(\overline S_n < \delta n w_n - n w_{\lfloor
n^\varepsilon
\rfloor}\bigr) %
\leq\mathbb{P}(\overline S_n < - \delta n w_n).
\end{equation}
By construction, $w_n$ is slowly varying and $a_n = L(n) n$ where
$L(n)$ is
slowly varying. Hence $| V_i | \leq\max(w_{\lfloor n^\varepsilon
\rfloor},
b_n) = L(n) n^{\varepsilon}$ where $L(n)$ is another slowly varying sequence.
By the Hoeffding inequality, we get from (\ref{eq:Sn})
\[
\mathbb{P}(\overline S_n < - \delta n w_n) %
\leq\exp%
\biggl( - \frac{ \delta^2 n ^2 w_n ^2 } { n L (n) ^2 n ^{2 \varepsilon
}} \biggr) %
= \exp( - \wt L(n) n ^{1 - 2 \varepsilon} ),
\]
where $\wt L(n)$ is a slowly varying sequence. Since $\varepsilon<
1/2$ we
obtain (\ref{eq:bcSn}) and thus (\ref{eq:llninfa1}).


\section{Properties of the limiting spectral distributions}
\label{sec:spec}

Recall that $\mu_\alpha$ is characterized by the Cauchy--Stieltjes transform
$m_{\mu_\alpha}(z) = \mathbb{E} h(z)$, $z\in\mathbb{C}_+$, where
$h(z)$ is the
random variable $h(z)=\langle\delta_\eset, (\bT-zI)^{-1}\delta
_\eset\rangle$ [cf. (\ref{meancon})]. The main novelty in our analysis
of the LSD $\mu
_\alpha$ with
respect to previous works \cite{benarous-guionnet,belinschi} is that
we can
work here with the distribution of $h(z)$ rather than only with its
expectation.

\subsection{Recursive distributional equation}\label{rdesub}
%
The symbol $ \stackrel{d}{=}$ stands for
equality in distribution. The following result is at the heart of our
analysis of the LSD~$\mu_\alpha$.
\begin{theorem}[(Recursive distributional equation)] \label{th:BC} For
all $z \in
\mathbb{C}_+$, the random variable
\[
h(z)=\langle\delta_\eset, (\bT-zI)^{-1}\delta_\eset\rangle
\]
satisfies to $h(-\bar z) = -\bar h (z)$ and
%
%
\begin{equation} \label{eq:RDE} %
h(z) \stackrel{d}{=} - \biggl( z + \sum_{k \in\mathbb{N}} \xi_k
h_k (z)
\biggr)^{-1},
\end{equation}
where $(h_k)_{k \in\mathbb{N}}(z)$ are i.i.d. with the same law of
$h(z)$, and $\{\xi_k\}_{k \in\mathbb{N}}$ is an independent Poisson point
process with intensity $\frac{\alpha}{2} x ^{-\alpha/2 -1}\,dx$ on
$(0,\infty)$.
\end{theorem}
%
%
\begin{pf}
Since the PWIT is bipartite, the property $h(-\bar z) = -\bar h (z)$ is
a consequence of Lemma~\ref{le:bipartite}. We are left with the RDE
(\ref{eq:RDE}). This can be interpreted as an operator version of the Schur
complement formula (see, e.g., Proposition 2.1 in Klein~\cite{klein}
for a
similar argument). Denote, as usual, by $k \in\mathbb{N}$ the
descendants of
the root $\eset$, and let $\cT^{(k)}$ denote the subtree rooted at
$k$ (the set
of vertices of $\cT^{(k)}$ is then $k \mathbb{N}^f$). We have the
direct sum
decomposition $\mathbb{N}^f = \{\eset\} \cup\bigcup_{k} k \mathbb
{N}^f$. We
define $\bT^{(k)}$ as the projection of $\bT$ on $k \mathbb{N}^f$. Its
skeleton is thus $\cT^{(k)}$. Finally, define the operator $\bU$ on
$\cD$ by
its matrix elements
\[
u_{k}:= \langle\delta_\eset, \bU\delta_k \rangle=
\langle\delta_k, \bU\delta_\eset\rangle= \langle\delta_\eset,
\bT\delta_k
\rangle\vadjust{\goodbreak}
\]
for all $k \in\mathbb{N}$ (offsprings of $\eset$) and $\langle
\delta_\mathbf{u}, \bU\delta_\bv\rangle=0$ otherwise. In this way
we have
\[
\bT= \bU+ \tbT\qquad\mbox{with } \tbT= \bigoplus_{k \in
\mathbb{N}}
\bT^{(k)}.
\]
As $\bT$, each $\bT^{(k)}$ can be extended to a self-adjoint operator,
which we denote again by $\bT^{(k)}$. Therefore\vspace*{1pt} $\tbT$ is
self-adjoint. We shall write $R(z)=(\bT- z I ) ^{-1}$ and $\wt R(z) =
(\tbT- z I ) ^{-1} $ for the associated resolvents, $z\in\mathbb{C}_+$.
These operators satisfy the resolvent identity
%
%
\begin{equation} \label{resid}
\wt R(z) (\bT-\tbT) R(z) = \wt R(z) - R(z) .
\end{equation}
Set $\wt R_{\mathbf{u},\bv}(z):=\langle\delta_\mathbf{u},\wt R(z)
\delta_\bv\rangle$ and
$R_{\mathbf{u},\bv}(z):=\langle\delta_\mathbf{u},R(z) \delta_\bv
\rangle$. Observe that\break $\wt R_{\eset,\eset}(z) = - z^{-1}$ and that the
direct sum decomposition\vspace*{1pt} $\mathbb{N}^f = \{\eset\}
\cup\bigcup_{k} k \mathbb{N}^f$ implies $\wt R_{k,l}(z)= 0$ for $k\neq
l$. Similarly we have that $\wt R_{\eset ,k}(z) = 0 = \wt
R_{k,\eset}(z)$ for every $k\in\mathbb{N}$. From (\ref{resid}) we then
obtain, for $k\in\mathbb{N}$,
\[
\wt R_{k,k}(z) u_k R_{\eset,\eset}(z) = - R_{k,\eset}(z) .
\]
It follows that
\begin{eqnarray*}
\langle\delta_\eset,\wt R(z) (\bT-\tbT) R(z)\delta_\eset
\rangle
&=& \sum_{k\in\mathbb{N}}\wt R_{\eset,\eset}(z) u_k R_{k,\eset}(z)\\
&=& -
\sum_{k\in\mathbb{N}} \wt R_{\eset,\eset}(z)\wt R_{k,k}(z) u_k^2
R_{\eset,\eset}(z) .
\end{eqnarray*}
From (\ref{resid}) we then conclude that
\[
R_{\eset,\eset}(z) = \frac{\wt R_{\eset,\eset}(z)}{1- \wt
R_{\eset,\eset}(z)\sum_{k\in\mathbb{N}} \wt R_{k,k}(z) u_k^2} .
\]
Or, using $\wt R_{\eset,\eset}(z) = -z^{-1}$,
\[
R_{\eset,\eset}(z) = -\biggl(z +
\sum_{k\in\mathbb{N}} \wt R_{k,k}(z) u_k^2\biggr) ^{-1} .
\]
Then (\ref{eq:RDE}) follows from the recursive construction of the PWIT:
$\cT^{(k)}$ are i.i.d. with distribution $\cT$ and therefore $\wt
R_{k,k}(z)$ are i.i.d. with the same law of $R_{\eset,\eset}(z)$, for
every $z\in\mathbb{C}_+$.
\end{pf}

Concerning the uniqueness of the solution to the RDE (\ref{eq:RDE}) we can
establish the following useful result. For $z = i t $, with $ t > 0 $, the
identity, $h(-\bar z) = -\bar h (z)$ reads $\Re h (it) = 0$. Thus, the
equation satisfied by $g(it) = \Im h(it) \geq0$ is
%
%
\begin{equation} \label{eq:RDEg} %
g(it) %
\stackrel{d}{=} \biggl( t + \sum_{k \in\mathbb{N}} \xi_k g_k (it)
\biggr)
^{-1}.
\end{equation}
\begin{lem}[(Uniqueness of solution for the RDE)]\label{le:itunique}
For each $t>0$, there exists a unique probability measure $L^{it}$ on
$\mathbb{R}_+$, solution of (\ref{eq:RDEg}).
\end{lem}
\begin{pf}
Set $\beta=\alpha/2$. If $(Y_k)$ is an i.i.d. sequence of nonnegative
random variables, independent of $\{ \xi_k\}_{k \in\mathbb{N}}$,
such that
$\mathbb{E} [ Y_1 ^\beta] < \infty$ then it is well known that
\[
\sum_k
\xi_k Y_k \stackrel{d}{=} \sum_k \xi_k (\mathbb{E} [ Y_1^\beta
])^{1/\beta}
\]
(see, e.g.,~\cite{talagrand2003}, Lemma 6.5.1, or (\ref{eq:laplace})
below). This implies the unicity for (\ref{eq:RDEg}) provided that
the equation satisfied by $ \mathbb{E}[ g(it)^\beta]$ has a unique solution.
Recall the formulas of Laplace transforms, for $y \geq0$, $\eta>0$
and $0
< \eta< 1$, respectively,
%
%
\begin{eqnarray} \label{eq:gammaLaplace} %
y^{-\eta} %
&=& \Gamma(\eta)^{-1} \int_0 ^\infty x ^{\eta-1} e^{- x y} \,dx\quad
\mbox{and} \nonumber\\[-8pt]\\[-8pt]
y^{\eta} %
&=& \Gamma(1-\eta)^{-1} \eta\int_0 ^\infty x ^{-\eta-1} (1 - e^{- x y})
\,dx.\nonumber
\end{eqnarray}
From the L\'{e}vy--Khinchine formula we deduce that, with $s \geq0$,
%
%
\begin{eqnarray} \label{eq:laplace}
\mathbb{E} \exp\biggl( - s \sum_{k} \xi_k Y_k \biggr)
& = & \exp\biggl( \mathbb{E} \int_0 ^\infty( e^{-x s Y_1} - 1 ) \beta
x ^{-\beta- 1} \,dx \biggr) \nonumber\\[-8pt]\\[-8pt]
&=& \exp\bigl( - \Gamma(1-\beta) s^\beta\mathbb{E}
[Y_1^\beta]\bigr).\nonumber
\end{eqnarray}
From (\ref{eq:RDEg}), $\mathbb{E}[ g(it)^\beta]$ is
the solution of the
equation in $y$:
\[
y = \frac{1}{\Gamma(\beta)} \int_0 ^\infty x ^{\beta-1} e^{-tx}
e^{- x ^{\beta} \Gamma(1-\beta) y} \,dx.
\]
The last equation has a unique solution for any $t\geq0$. Indeed, the
function from $\mathbb{R}_+$ to $\mathbb{R}_+$
\[
\varphi\dvtx y \mapsto\frac{1}{\Gamma(\beta)} \int_0 ^\infty x
^{\beta-1} e^{-tx} e^{- x ^{\beta} \Gamma(1-\beta) y}\,dx
\]
tends
to $0$ as $y\to\infty$, and it is decreasing since
\[
\varphi'(y) %
= - \frac{\Gamma(1- \beta) }{\Gamma(\beta)}\int_0 ^\infty x
^{2\beta-1}
e^{-tx} e^{- x ^{\beta} \Gamma(1-\beta) y}\,dx.
\]
Thus $\varphi$ has a unique fixed point.
\end{pf}

Before going into the proof of Theorem~\ref{th:mua}, we introduce some
notation. Let $\beta= \alpha/2$ as above, and let $\cK_\alpha$
denote the set
of probability measures on $(0,\infty)$ with finite $\beta$ moment.
We define
the map $\Psi$ on probability measures on $\mathbb{R}_+ \cup\{\infty
\}$,
where $\Psi(Q)$ is the law of
%
%
\begin{equation}\label{eq:Psi}
Z = \biggl( \sum_{k \in\mathbb{N}} \xi_k Y_k \biggr)^{-1},
\end{equation}
with $(Y_k,k \in\mathbb{N})$ i.i.d. with law $Q$ independent of $\Xi=
\{\xi_k\}_{k \in\mathbb{N}}$ a Poisson point process on $\mathbb
{R}_+$ of
intensity $\beta x ^{-\beta- 1}\,dx$.
\begin{lem}\label{le:Psi}
$\Psi$ satisfies the following:
\begin{enumerate}[(iii)]
\item[(i)] $\Psi$ is a map\vspace*{1pt} from $\cK_\alpha$ to $\cK_\alpha$. Let
$(P_n)_{n
\in\bN}$ and $P$ in $\cK_\alpha$, if $\lim_{n\to\infty} \int
x^\beta\,
dP_n = \int x^\beta \,dP$ then $\Psi(P_n)$ converges weakly to $\Psi(P)$
and $\lim_{n\to\infty} \int x^\beta \,d\Psi(P_n) = \int x^\beta \,d\Psi(P)$.
\item[(ii)] The unique fixed point of $\Psi$ in $\cK_\alpha$ is the
law of
$1/S$ where $S$ is the one-sided $\beta$-stable law with Laplace transform
$ \mathbb{E} \exp( - t S ) = \exp( -t^{\beta}\times\sqrt{ \Gamma
(1 +
\beta)/\Gamma(1 -\beta)} ) $, $t\geq0$.
\item[(iii)] $\mathbb{E} S^{-\beta} = (\Gamma(\beta+1)\Gamma( 1 -
\beta) )
^{-1/2}$.
\end{enumerate}
\end{lem}
\begin{pf}
As in the proof of Lemma~\ref{le:itunique}, we get
\begin{eqnarray*}
\mathbb{E} Z^\beta& = & \mathbb{E} \biggl( \sum_k \xi_k Y_k
\biggr)^{-\beta} \\
& = & \mathbb{E} \frac{1}{\Gamma(\beta)} \int_0 ^\infty x^{\beta-
1} e^{-x \sum_k \xi_k Y_k} \,dx \\
& = & \frac{1}{\Gamma(\beta)} \int_0 ^\infty x^{\beta- 1} e^{-
x^\beta\Gamma( 1 - \beta) \mathbb{E} Y_1^\beta} \,dx \\
& = & \frac{1}{\beta\Gamma(\beta)} \int_0 ^\infty e^{- s \Gamma(
1 - \beta) \mathbb{E} Y_1^\beta} \,ds\\
& = & \bigl(\Gamma(\beta+1)\Gamma( 1 - \beta) \mathbb{E} Y_1^\beta
\bigr)^{-1} ,
\end{eqnarray*}
[in the last line we have used the identity $z \Gamma(z) = \Gamma(z+1)$].
Therefore, $\Psi$ is a map from $\cK_\alpha$ to $\cK_\alpha$. Also
as a
consequence of (\ref{eq:laplace})
\[
\mathbb{E} \exp( -t Z^{-1} )%
= \exp\bigl( - t^\beta\Gamma( 1 - \beta) \mathbb{E} Y_1^\beta\bigr).
\]
Statement (i) follows from the continuity of the map $x \mapsto1/x$ in
$(0,\infty)$. If $Z$ is a fixed point of $\Psi$ then from the computation
above $\mathbb{E} Z^\beta= (\Gamma(\beta+1)\Gamma( 1 - \beta) )
^{-1/2}$.
Finally, from (\ref{eq:laplace}) we obtain for all $t\geq0$,
\[
\mathbb{E} \exp( -t Z^{-1} ) %
= \exp\bigl( - t^\beta\Gamma( 1 - \beta) \mathbb{E} Z^\beta\bigr) %
= \exp\Biggl( - t^{\beta} \sqrt{\frac{ \Gamma(1 + \beta)}{ \Gamma
(1 -\beta)
} } \Biggr) .
\]
\upqed\end{pf}

\subsection{\texorpdfstring{Proof of Theorem \protect\ref{th:mua}\textup{(i)}}
{Proof of Theorem 1.6(i)}}

From Theorem~\ref{th:BC}, for $z \in\mathbb{C}_+$,
\[
m_{\mu_\alpha} (z) = \mathbb{E} h (z),
\]
where $h$ solves RDE (\ref{eq:RDE}). Set $f(z) = \Re h(z)$ and $g(z) =
\Im
h(z)$. For $z = u + i v \in\mathbb{C}_+$, $f$ and $g$ satisfy the RDE
\[
f(z) \stackrel{d}{=} %
- \frac{ u + \sum_k \xi_k f_k (z) } { ( u + \sum_k \xi_k
f_k(z) )^2 + ( v + \sum_k \xi_k g_k(z) )^2}
\]
and
\[
g(z) \stackrel{d}{=} \frac{ v + \sum_k \xi_k g_k (z) } { ( u
+ \sum_k \xi_k f_k(z) )^2 + ( v + \sum_k \xi_k g_k
(z))^2} .
\]
By construction, $0 \leq g(z) \leq1/v$, thus the law of $g(z)$ is in
$\cK_\alpha$. If the stochastic domination of $P$ by $Q$ is denoted
by $P
\leq_{st} Q$, we have
%
%
\begin{equation}\label{eq:gborne}
g(z) \leq_{st} \biggl(v+ \sum_k \xi_k g_k (z)\biggr)^{-1} %
\leq_{st} \biggl( \sum_k \xi_k g_k (z)\biggr)^{-1}.
\end{equation}
[In fact, we also have $| h(z) | \leq_{st} (\sum_k \xi_k g_k
(z))^{-1}$.] Using\vspace*{1pt} the computation in Lem\-ma~\ref{le:Psi}, we obtain
$\mathbb{E} g(z) ^\beta\leq(\Gamma(\beta+1)\Gamma( 1 - \beta)
\mathbb{E} g(z) ^\beta
)^{-1}$. Thus
%
%
\begin{equation}\label{eq:tension}
\mathbb{E} g(z) ^\beta\leq\frac{1}{\sqrt{\Gamma(\beta+1)\Gamma(
1 - \beta) ) }}.
\end{equation}
Again, the formula $y^{-\eta} = \Gamma(\eta)^{-1} \int_0 ^\infty x
^{\eta
-1} e^{- x y} \,dx$, for $y \geq0$, $\eta>0$, gives
%
%
\begin{equation}\label{eq:gmoment}
\mathbb{E} \biggl[\biggl( \sum_k \xi_k g_k (z) \biggr)^{-\eta
}\biggr] =
\frac{1}{\Gamma(\eta)}
\int_0 ^\infty x ^{\eta-1} e^{- x^\beta\Gamma(1-\beta) \mathbb{E}
g(z) ^\beta} \,dx .
\end{equation}
We now study the weak limit of $g(u + iv)$ when $v \downarrow0$, $u
\in
\mathbb{R}$. Equation (\ref{eq:tension}) implies tightness, so let
$g(u +i 0)$ be
a weak limit. If this limit is nonzero then $\mathbb{E} g^\beta(u +
i0) > 0$,
and equations (\ref{eq:gborne})--(\ref{eq:gmoment}) imply for all
$\eta>0$
and $u \in\mathbb{R}$,
\[
\limsup_{u + i v \dvtx v \downarrow0} \mathbb{E} g^\eta(u + i v) <
\infty.
\]
Since $\mathbb{E} h(z)$ is the Cauchy--Stieltjes transform of $\mu
_\alpha$,
taking $\eta= 1$, we deduce that $\mu_\alpha$ is absolutely
continuous (see, e.g.,~\cite{simon05}, Theorem 11.6).

\subsection{\texorpdfstring{Proof of Theorem \protect\ref{th:mua}\textup{(ii)}}
{Proof of Theorem 1.6(ii)}}

In view of~\cite{simon05}, Theorem 11.6, it is sufficient to show that
%
%
\begin{equation}\label{eq:limg0}
\lim_{t \downarrow0} \mathbb{E} g(it) %
= \Gamma\biggl(1+ \frac{1}{\beta}\biggr) \biggl( \frac{ \Gamma(1
+ \beta)}{ \Gamma(1 - \beta) } \biggr)^{{1}/({2 \beta})}.
\end{equation}
As above, (\ref{eq:tension}) implies the tightness of $(g(it), t >
0)$. So
let $g(i 0)$ be a weak limit. It is in $\cK_\alpha$ and, by continuity,
$g(i0)$ is solution of the RDE
\[
g(i0) \stackrel{d}{=} \biggl( \sum_k \xi_k g_k (i0)\biggr)^{-1}.\vadjust{\goodbreak}
\]
By Lemma~\ref{le:Psi}, $g(i0) \stackrel{d}{=} 1/ S$, and (\ref{eq:gmoment})
gives
\[
\mathbb{E} g(i0) = \int_0 ^\infty e^{- x^\beta\sqrt{{\Gamma
(1-\beta)
}/{\Gamma(1+\beta)}} } \,dx = \frac{1}{\beta}
\Gamma\biggl(\frac{1}{\beta}\biggr) \biggl( \frac{ \Gamma(1 +
\beta)}{
\Gamma(1 - \beta) } \biggr)^{{1}/({2 \beta})}.
\]
Using the identity $z \Gamma(z) = \Gamma(z+1)$, we get (\ref{eq:limg0}).

\subsection{\texorpdfstring{Proof of Theorem \protect\ref{th:mua}\textup{(iii)}}
{Proof of Theorem 1.6(iii)}}

We start with a Tauberian-type theorem for the Cauchy--Stieltjes
transform of
symmetric probability measures. As usual, let $m_\mu$ denote the
Cauchy--Stieltjes
transform of a symmetric probability measure $\mu$ on~$\mathbb{R}$.
Then, for
all $t > 0$, $m_\mu(it ) \in i \mathbb{R}_+$ and
\[
\Im m_\mu(it) =
\int_{-\infty}^\infty\frac{ t }{ t^2 + x^2} \mu(dx)= 2 \int_0
^\infty\frac{
t }{ t^2 + x^2} \mu(dx).
\]

\begin{lem}[(Tauberian-like lemma)]\label{le:taub}
If $L$ is slowly varying and $0 < \alpha< 2$, the following are equivalent:
as $t$ goes to $+\infty$
%
%
\begin{eqnarray}
\label{eq:taub1}
\mu( (t,\infty))
& \sim & L(t) t ^{- \alpha}, \\
\label{eq:taub2}
\Im m_\mu(it) - t ^{-1}
& \sim & - \Delta(\alpha) L(t) t ^{- \alpha-1}
\end{eqnarray}
with $ \Delta(\alpha) = 2\alpha\int_0 ^\infty\frac{ x^{1 - \alpha
} } { 1 + x ^2 } \,dx$.
\end{lem}
\begin{pf*}{Sketch of Proof of Lemma~\ref{le:taub}} The proof is an
adaptation of the proof of the Karamata's Tauberian theorem in
\cite{bingham}, pages 37 and 38. Let $\cM$ denote the set of symmetric
measures on $\mathbb{R}$ such that $\int_0^\infty\min(1, x^2) \mu
(dx) <
+ \infty$. On $\cM$, define the transform
\[
\cS\mu\dvtx t \mapsto\int_0 ^\infty\frac{2 x^ 2} { t^2 + x ^2 } \mu
( dx).
\]
Note that $\cS\mu(t) = 1 - t \Im m_\mu(it) = 1 + i t m_\mu(it)$. Recall
that the Cauchy--Stieltjes transform characterizes the measure. Thus if for
all $t > 0$, $(\cS\mu_n (t))_{n \in\mathbb{N}} $ converges to $\cS
\mu$,
then $(\mu_n)_{n \in\mathbb{N}} $ converges to $\mu$ over all bounded
continuous function with $0$ outside the support. Now, assume that
(\ref{eq:taub2}) holds, namely
%
%
\begin{equation}\label{scond} \cS\mu(t)
\sim\Delta(\alpha) L(t) t ^{- \alpha}.
\end{equation}
Since $\lim_{x \to\infty} L( t x ) / L(t) = 1$, we deduce that for
all $t
>0$, as $x \to\infty$
\[
\frac{ \cS\mu(x t) }{ L(x) x ^{-\alpha} } \to\Delta(\alpha) t
^{- \alpha}.
\]
The left-hand side is the $\cS$ transform of the measure $\mu_x (dy)
= \mu(x\,
dy ) /\break ( L(x) x ^{-\alpha} )$ while the right-hand side is the $\cS$
transform of $\mu_{\infty} (dy ) =\break \alpha|y | ^{- \alpha-1} \,dy$,
thus
\[
\frac{ \mu(( x, \infty) ) }{ L(x) x ^{-\alpha} } = \mu_x (
(1,\infty)
) \quad\to\quad\mu_\infty(1,\infty) = 1.
\]
We get precisely (\ref{eq:taub1}). The
reciprocal implication can be proved similarly (see~\cite{bingham},
pages 37 and 38) [it is straightforward for $L(t) = c$, the case that
we will
actually use].
\end{pf*}

We now come back to the RDE (\ref{eq:RDEg}) and define $Q(t) = \mathbb{E}[
g(it)^\beta]$. From (\ref{eq:RDEg}), we have a.s. $t g(it) \leq1$. Note
also, from a.s. $\sum_k \xi_k g_k (it) \leq t^{-1} \sum_k \xi_k $,
that a.s.
$\lim_{t \to+\infty} t g(it) = 1$. The dominated convergence theorem leads
to
%
%
\begin{equation}
\label{eq:asyG}
\lim_{t \to\infty} t ^{ \beta}
Q(t) =1.
\end{equation}
Moreover, as already pointed in Lemma~\ref{le:itunique},
\[
\sum_k \xi_k g_k
(it) \stackrel{d}{=} Q(t)^{1/\beta} \sum_k \xi_k.
\]
We deduce, with $C(t) = ( t
Q(t)^{1/\beta})^{-1/2}$, that
%
%
\begin{eqnarray} \label{eq:asyma}
\Im m_{\mu_\alpha} ( it) &=& \mathbb{E} g (it) = \mathbb{E} \frac
{ t } { t^2 + t Q(t)^{1/\beta} \sum_k \xi_k} \nonumber\\
& = & C(t) \mathbb{E} \frac{ t C(t)} { (t C(t))^2 + \sum_k \xi_k}
\\
& = & C(t) \Im m_{\cL(Y)} ( i C(t) t ),\nonumber
\end{eqnarray}
where $\cL(Y)$ is the law of
\[
Y = \varepsilon\sqrt{\sum_k \xi_k},
\]
and
$\varepsilon$ is independent of $\{\xi_k\}_k$, $\mathbb{P}
(\varepsilon=1 ) = \mathbb{P}
(\varepsilon= - 1 ) = 1/2$. We have
\[
\mathbb{P} ( Y > t ) = \frac1 2 \mathbb{P} \biggl( \sum_k \xi_k >
t^2\biggr).
\]
By (\ref{eq:laplace}), as $s \downarrow0$, $\mathbb{E} \exp( -s
\sum_k \xi_k
) = \exp( - s^\beta\Gamma(1 - \beta) ) \sim1 - s^\beta\Gamma(1 -
\beta)$. Using~\cite{bingham}, Corollary 8.7.1, we obtain $\mathbb
{P} ( \sum_k
\xi_k > t ) \sim t^{-\beta}$ and
\[
\mathbb{P} ( Y > t ) \sim\frac{t^{-\alpha}}{2}.
\]
By Lemma~\ref{le:taub}, $ \Im m_{\cL(Y)} ( i t ) - t ^{-1} \sim-
\frac{t^{-\alpha-1}}{2} \Delta(\alpha)$. Thus by
(\ref{eq:asyG}) and (\ref{eq:asyma}), 
\[
\Im m_{\mu_\alpha} ( it) - t ^{-1} \sim- \frac{t^{-\alpha
-1}}{2}\Delta(\alpha).
\]
Theorem~\ref{th:mua}(iii) now follows from Lemma~\ref{le:taub}.
\begin{rem}
In the proof of Lemma~\ref{le:itunique}, we have seen that the distribution
of $g(it) = \Im h(it)$ was function\vadjust{\goodbreak} of $Q(t) = \mathbb{E} [ g^\beta
(it) ]$ which
satisfies the equation
\[
Q(t) = \frac{1}{\Gamma(\beta)} \int_0 ^\infty x ^{\beta-1} e^{-tx}
e^{- x
^{\beta} \Gamma(1-\beta) Q(t)} \,dx = f_\beta( t, Q(t)).
\]
We could push further our investigation at $t =0$ and compute the derivative
of $Q$ at\vspace*{1pt} $t =0$: $Q'(0) = - f_{\beta+1} (0,Q(0)) - \Gamma( 1 -
\beta)f_{2\beta} (0,Q(0)) Q'(0)$, with $Q(0) = (\Gamma(\beta
+1)\Gamma( 1 -
\beta) )  ^{-1/2}$. There should be no obstacle for computing by recursion
the successive derivatives of $Q(t)$ at $t=0$. We would then obtain a series
expansion of the partition function $\mu_\alpha((-\infty,t))$ in a
neighborhood of $0$.
\end{rem}

\subsection{\texorpdfstring{Proof of Theorem \protect\ref{th:muabis}: $\tilde{\mu}_\alpha$, $\alpha\in(0,1)$}
{Proof of Theorem 1.7: mu alpha, alpha in (0,1)}}

As in (\ref{gammal}), let $\mathbf{p}_{\ell}$ denote the return
probability after
$\ell$ steps starting from the root $\eset$, for the random walk on
the PWIT
with transition kernel $\bK$ given by (\ref{kappone}). In particular,
$\gamma_{\ell} = \mathbb{E} \mathbf{p}_{\ell}$ is the $\ell$th moment
of the LSD
$\wt\mu_\alpha$.
\begin{pf*}{Proof of Theorem~\ref{th:muabis}\textup{(i)}}
For the first part, we shall show that there exists $\delta>0$ such
that for
any $\varepsilon\in(0,1/2]$ and any $n$
%
%
\begin{equation}\label{hm2}
\gamma_{2n}\geq\delta\varepsilon^\alpha(1-\varepsilon)^{2n} .
\end{equation}
Theorem~\ref{th:muabis} (i) follows by choosing $\varepsilon=1/2n$.
To prove (\ref{hm2}) we use the simple bound $\mathbf{p}_{2n}\geq
(\bK(\eset,1)\bK(1,\eset))^n$, which states that to
come back to
the root in $2n$ steps the walk can move to the child with the highest
weight, with probability $\bK(\eset,1)$, go back to the root, with
probability $\bK(1,\eset)$, and repeat this $n$ times. Taking expectation,
it follows that
%
%
\begin{equation}\label{hm3}
\gamma_{2n} \geq\mathbb{E}[(\bK(\eset,1)\bK(1,\eset
))^n] .
\end{equation}
Therefore (\ref{hm2}) holds if the event
\[
A_\varepsilon= \{\bK(\eset,1)\geq(1-\varepsilon) \mbox{ and }
\bK(1,\eset)\geq(1-\varepsilon)\}
\]
has probability at least $\delta\varepsilon^\alpha$, for some
$\delta>0$
and for any $\varepsilon\in(0,1/2]$.

Let $(x_i)_i$ denote the realization of the PPP at the root $\eset$,
that is,
$x_1>x_2>\cdots$ are the points of a PPP on $(0,\infty)$ with intensity
measure $\alpha x^{-\alpha-1}\,dx$. We set $\phi:=\sum_{i=1}^\infty
x_i$ and let
$\phi'$ denote an independent copy of $\phi$. We can use the representation
$\bK(\eset,1) = x_1/\phi$ and $\bK(1,\eset)=x_1/(x_1+\phi')$. Therefore,
\begin{eqnarray*}
\mathbb{P}(A_\varepsilon) & = &
\mathbb{P}\bigl(x_1\geq(1-\varepsilon)\phi, x_1\geq
(1-\varepsilon)(x_1+\phi')\bigr)\\
& = & \mathbb{P}\biggl(x_1\geq(1-\varepsilon)\phi, \phi'\leq
\frac{\varepsilon x_1}{(1-\varepsilon)}\biggr)\\
&\geq&\mathbb{P}\bigl(x_1\geq(1-\varepsilon)\phi, x_1\geq
\varepsilon^{-1} , \phi'\leq1\bigr) .
\end{eqnarray*}
Let $\delta_1:= \mathbb{P}(\phi\leq1) = \int_0^1 f(t) \,dt>0$,
where $f(t)$
denotes the density of $\phi$. The function $f(t)$ can be obtained
from its
Laplace transform, which is given by the known identity $\mathbb
{E}[e^{-u\phi}]
=e^{-\Gamma(1-\alpha)u^{\alpha}}$, $u>0$ (see~\cite{MR1434129}, Proposition
10, or (\ref{eq:laplace}) with $\beta$ replaced by
$\alpha$ and
$Y_k = 1$). Since $\phi'$ is independent of $(x_i)$ we obtain
\[
\mathbb{P}(A_\varepsilon) %
\geq\delta_1 \mathbb{P}\bigl(x_1\geq(1-\varepsilon)\phi,
x_1\geq
\varepsilon^{-1}\bigr) .
\]
To estimate the last quantity we observe that if $\wt x$ is a size-biased
pick from $(x_i)$, then $x_1\geq\wt x$. We recall that $\wt x$ is a random
variable such that, given the sequence $(x_i)$ the probability that
$\wt x$
equals $x_i$ is $x_i/\phi$. It is not hard to check (see, e.g., \cite
{PerPitYor},
Lemma 2.2) that the random variable $\wt x$ has a probability density
on $(0,\infty)$ given by
%
%
\begin{equation}\label{sz}
\alpha x^{-\alpha-1}\int_0^\infty f(t) \frac{x}{x+t} \,dt ,
\end{equation}
where $f(t)$ is the density of the variable $\phi$. Therefore,
\begin{eqnarray*}
&&\mathbb{P}\bigl(x_1\geq(1-\varepsilon)\phi, x_1\geq
\varepsilon^{-1}\bigr)\\
&&\qquad\geq\mathbb{P}\bigl(\wt x\geq(1-\varepsilon)\phi, \wt x\geq
\varepsilon^{-1}\bigr)\\
&&\qquad=\alpha\int_0^\infty dt\, f(t)\int_0^\infty
\,dx\, x^{-\alpha-1} \frac{x}{x+t}
\ind_{\{x\geq(1-\varepsilon)(x+t)\}} \ind_{\{x\geq\varepsilon
^{-1}\}} \\
&&\qquad\geq \alpha\int_0^1 dt\, f(t)\int_0^\infty
\,dx\, x^{-\alpha-1} (1-\varepsilon) \ind_{\{x\geq\varepsilon^{-1}\}
}\\
&&\qquad = \delta_1 (1-\varepsilon) \varepsilon^\alpha.
\end{eqnarray*}
In conclusion, $\mathbb{P}(A_\varepsilon)\geq\delta_1^2
(1-\varepsilon) \varepsilon^\alpha\geq
\frac12 \delta_1^2 \varepsilon^\alpha$, and the claim (\ref
{hm2}) follows.

It remains to show that $\liminf_{\alpha\nearrow1}\gamma_2>0$. If $(x_i)$,
$\wt x$, and $\phi$ are as above and if $\phi'$ is independent of the
sequence $(x_i)$ and identical in law to the random variable~$\phi$, then
\[
\gamma_2
=\mathbb{E}\biggl[\sum_{i}\frac{x_i}{\phi}\frac{x_i}{x_i+\phi
'}\biggr]
=\mathbb{E}\biggl[\frac{\wt x}{\wt x+\phi'}\biggr]
=\int_0^\infty\alpha
x^{1-\alpha}\biggl(\int_0^\infty\frac{f(t)}{x+t} \,dt\biggr)^2 \,dx.
\]
Now, from the Laplace transform $\mathbb{E}[e^{-u\phi}]
=e^{-\Gamma(1-\alpha)u^{\alpha}}$ we have the identity
\[
\int_0^\infty\frac{f(t)}{x+t} \,dt %
=\int_0^\infty e^{-\Gamma(1-\alpha)u^\alpha-ux} \,du .
\]
This gives
\begin{eqnarray*}
\gamma_2
&=&\alpha\Gamma(2-\alpha)\int_0^\infty\int_0^\infty
e^{-\Gamma(1-\alpha)(u^\alpha+v^\alpha)} (u+v)^{-2+\alpha} \,du\,
dv\\
&=&\frac{\alpha\Gamma(2-\alpha)}{\Gamma(1-\alpha)}
\int_0^\infty\int_0^\infty e^{-t^\alpha-s^\alpha}
(t+s)^{-2+\alpha} \,ds \,dt .
\end{eqnarray*}
Finally, the desired result follows from the bounds (for absolute
constants $c_1$, \mbox{$c_2>0$})
\begin{eqnarray*}
\int_0^\infty\int_0^\infty e^{-t^\alpha-s^\alpha}
(t+s)^{-2+\alpha} \,ds \,dt %
&\geq& e^{-2} \int_0^1 \int_0^1 (t+s)^{-2+\alpha} \,ds\,dt \\
&\geq&\frac{c_1}{1-\alpha}
\end{eqnarray*}
and
\[
\Gamma(1-\alpha) =\int_0^\infty t^{-\alpha}e^{-t} \,dt \leq
\int_0^1 t^{-\alpha} \,dt+\int_1^\infty e^{-t} \,dt\leq
\frac{c_2}{1-\alpha}.
\]
\upqed\end{pf*}
\begin{pf*}{Proof of Theorem~\ref{th:muabis}\textup{(ii)}}
It is convenient to make here the dependence over $\alpha$ explicit in all
the notation. In particular, for every $\alpha\in(0,1)$, we denote by
$\bS_\alpha$ the operator $\bS$ given by (\ref{opsym}). These
operators are
defined on a common probability space, and are self-adjoint in $L^2 ( V)$.
Moreover, it follows from Section~\ref{sec:k01} that $\widetilde\mu
_\alpha=
\mathbb{E} \mu_{\alpha, \eset}$, where $\mu_{\alpha,\eset}$ is
the spectral
measure of $\bS_\alpha$ at the vector $\delta_{\eset}$. By the dominated
convergence theorem, in order to prove that $\alpha\mapsto\widetilde
\mu_\alpha$ is continuous in $(0,1)$, it is sufficient to show that
a.s. $\alpha\mapsto\mu_{\alpha,\eset}$ is continuous. From~\cite{reedsimon},
Theorem VIII.25(a),
it is in turn sufficient to prove that
for all $\bv\in V$, $\alpha\mapsto\bS_{\alpha} \delta_\bv$ is a
continuous map from $(0,1)$ to $L^2 ( V)$. From (\ref{opsym}), for all
$\mathbf{u}
\in V$, the map $\alpha\mapsto\bS_\alpha( \mathbf{u}, \bv) $ is
continuous. It
thus remains to check the uniform square integrability of
$(\bS_{\alpha} ( \bv, \mathbf{u}) )_{\mathbf{u}\in
V}$. We start with the
upper bound
\[
(\bS_{\alpha} ( \bv, \bv k ) ) ^2 = \frac{y_{\bv k} ^ {-1 /
\alpha} }{\rho_\alpha(\bv)} \frac{y_{\bv k} ^ {-1 / \alpha}
}{\rho_\alpha(\bv k )} \leq\frac{y_{\bv k}^{-1 / \alpha} }{\rho
_\alpha
(\bv)}.
\]
Then, notice that for all $\alpha\in(0, 1 - \varepsilon)$, one has
$y_{\bv k} ^ {-1 / \alpha} \leq\max(1, y_{\bv k} ^ {-1 /
(1-\varepsilon)})$ and $ \rho_\alpha(\bv) \geq\min( 1, y_{\bv1} ^{- 1/
( 1 - \varepsilon) })$. We may conclude
by recalling that a.s. $\lim_k y_{\bv k} / k =1$ and $y_{\bv1} >0$.
\end{pf*}
\begin{pf*}{Proof of Theorem~\ref{th:muabis}\textup{(iii)}}
As in the proof of Theorem~\ref{th:muabis}(ii), we make here the dependence
over $\alpha$ explicit in all the notation. It follows from Section~\ref{sec:k01}
\[
\int x^{2 \ell} \widetilde\mu_\alpha(dx) %
= \mathbb{E} \int x^{2 \ell} \mu_{\alpha, \eset} (dx) %
= \mathbb{E} \mathbf{p}_{\alpha,2 \ell} ,
\]
where the expectation is over the randomness of the PWIT. We introduce for
$\bv\in V$,
\[
V_\alpha( \bv) = \biggl( \frac{y_{\bv1} ^{-1/\alpha}}{ \sum_{k
\geq1}
y_{\bv k} ^{-1/\alpha}} , \frac{y_{\bv2} ^{-1/\alpha}}{ \sum_{k
\geq1}
y_{\bv k} ^{-1/\alpha}}, \ldots\biggr).
\]
By construction $ V_\alpha( \bv)$ is a PD($\alpha,0$) random variable.
Thus, by~\cite{MR1434129}, Corollary 18, as $\alpha\downarrow0$,
$V_\alpha
( \bv)$ converge weakly to the deterministic vector $(1,0,\ldots)$.
We may
thus write
\[
\bK_\alpha(1,\eset) = \frac{y_1^{-1/\alpha}}{ y_1^{-1/\alpha} +
y_{11}^{-1/\alpha}(1+\varepsilon_\alpha)},
\]
where as $\alpha$ goes to $0$, $\varepsilon_\alpha$ goes in
probability to
$0$. We define $U = \ind_{\{ y_{11} > y_{1} \}}$, so that $U$ is a symmetric
Bernoulli, that is, $\mathbb{P}(U = 0) = \mathbb{P}( U =1) = 1/2$. We
have proved
that in probability
\[
\lim_{\alpha\downarrow0} \bK_\alpha( \eset, 1 ) = 1
\quad\mbox{and}\quad
\lim_{\alpha\downarrow0} \bK_\alpha( 1, \eset) = U.
\]
In particular, 
\[
\lim_{\alpha\downarrow0} \int x^{2 \ell} \mu_{\alpha,\eset} (dx)
= U .
\]
Since $ \mu_{\alpha,\eset}$ is symmetric, 
\[
\lim_{\alpha\downarrow0} \mu_{\alpha,\eset} %
= \frac{U}{2} \delta_{-1} + (1-U) \delta_0 + \frac{U}{2} \delta_1.
\]
Taking expectation, we obtain the claimed statement on $ \widetilde\mu
_\alpha$.
\end{pf*}

\section{\texorpdfstring{Invariant measure: Proof of Theorem \protect\ref{th:inv}}
{Invariant measure: Proof of Theorem 1.8}}
\label{sec:inv}

We start with a lemma. Let $(X_1,\ldots,X_n)$, $X_1\geq\cdots\geq X_n$,
denote the ranked values of $\rho_1,\ldots,\rho_n$ and recall the
notion of
convergence in the space $\cA$, cf. Section~\ref{order}. We use the notation
$b_n:=a_{m_n}$, where $m_n=n(n+1)/2$.
\begin{lem}\label{ordst} For any $\alpha\in(0,2)$,
the sequence $b^{-1}_{n}(X_1,X_2,\ldots)$ converges in distribution to
$(x_1,x_1,x_2,x_2,\ldots)$, where $x_1>x_2>\cdots$ denote the ranked points
of the Poisson point process on $(0,\infty)$ with intensity
$\alpha x^{-\alpha-1}\,dx$.
\end{lem}
\begin{pf}
There are $m_n=n(n+1)/2$ edges, including self-loops. Let us denote by
$U_e$ the weight of edge $e\in\{1,\ldots,m_n\}$. The row sums are
given by
$\rho_i = \sum_{e\dvtx e\ni i} U_e$. We write $O_n$ for the set of off-diagonal
edges $e$, that is, edges of the form $e=\{i,j\}$ with $i\neq j$. Let
$U_{e_1}\geq U_{e_2}\geq\cdots$ denote the ranked values of the
i.i.d. random vector $(U_e)_{e\in O_n}$. Since there are $m_n-n$ edges in
$O_n$, an
application of Lemma~\ref{le:PoiExt}(i) yields convergence in distribution
%
%
\begin{equation}\label{sosh1}
b_n^{-1}(U_{e_1},U_{e_2},\ldots)
\mathop{\longrightarrow}^{d}_{n\to\infty}
(x_1,x_2,\ldots) .
\end{equation}
Each $e_i=\{u_i,v_i\}\in O_n$ identifies two row sums $\rho_{u_i}$ and
$\rho_{v_i}$. Set $\Delta_i =\break
\max\{\rho_{u_i}-U_{e_i},\rho_{v_i}-U_{e_i}\}$. Then, for every
$k\in\mathbb{N}$
and $\varepsilon>0$,
%
%
\begin{equation}\label{sosh}
\lim_{n\to\infty}
\mathbb{P}\Bigl(\max_{1\leq\ell\leq k}\Delta_\ell\geq
\varepsilon b_n\Bigr) = 0 .
\end{equation}
To prove this we use an estimate due to Soshnikov~\cite{MR2081462}. Let
$B_n$ denote the event that there exists no $i\in\{1,\ldots,n\}$ such that
\[
\Bigl\{\rho_i > b_n^{3/4 + \alpha/{8}} \mbox{ and } \rho_i
- \max_{j}
U_{i,j} > b_n^{3/4 + \alpha/{8}}\Bigr\}.
\]
Then, from~\cite{MR2081462} and~\cite{auffinger-benarous-peche}, Lemma 3,
one has
%
%
\begin{equation}\label{soshn}
\lim_{n\to\infty}\mathbb{P}(B_n)\to1 .
\end{equation}
Clearly, on the event $B_n$, if $\max_{1\leq\ell\leq k}\Delta_\ell
\geq
\varepsilon b_n$, then $U_{e_k}\leq b_n^{3/4 + \alpha/
{8}}$ which
has vanishing probability in the limit by (\ref{sosh1}). This proves
(\ref{sosh}).

For simplicity, we introduce the notation $R_{2\ell-1} =
\max\{\rho_{u_{\ell}},\rho_{v_{\ell}}\}$, $R_{2\ell} =
\min\{\rho_{u_{\ell}},\rho_{v_{\ell}}\}$. Therefore (\ref{sosh}) and
(\ref{sosh1}) prove that
%
%
\begin{equation}\label{sosh10}
b_n^{-1}(R_1,R_2,R_3,R_4,\ldots)
\mathop{\longrightarrow}^{d}_{n\to\infty}
(x_1,x_1,x_2,x_2,\ldots) .
\end{equation}
It remains to show that for every fixed $k$
%
%
\begin{equation}\label{star}
\lim_{n\to\infty} 
\mathbb{P}\biggl(\bigcup_{1\leq i\leq2k}\{R_i\neq X_i\}\biggr) = 0 .
\end{equation}
By construction, we have $X_i\geq R_i$ for $i=1,2$. On the event $B_n$
described above, to have $X_1>R_1$ or $X_2>R_2$ implies that there
exists an
edge $e\neq e_1$ such that $U_e \geq U_{e_1} - b_n^{3/4 +
\alpha/{8}}$. However, this event has vanishing probability by
(\ref{sosh1}) and the fact that $b_n^{\delta-1}\max_i U_{i,i}\to0$ in
probability for all sufficiently small $\delta>0$ (indeed by Lemma
\ref{le:PoiExt}, $a_n ^{-1} \max_i U_{i,i}$ converges weakly to the
Fr\'{e}chet distribution, see first comment after Lemma~\ref{le:PoiExt}).
Thanks to (\ref{soshn}) this shows that $\mathbb{P}(X_1>R_1 \mbox
{ or }
X_2>R_2) \to0$. Recursively, the probability of $X_{2i+1}> R_{2i+1}$ or
$X_{2i+2}> R_{2i+2}$ on the event $B_n\cap\{X_j=R_j, \forall
j=1,\ldots,2i\}$ vanishes as $n\to\infty$. Indeed, at each step we have
removed a row and a column corresponding to the largest off-diagonal weight
and we may repeat the same reasoning as above. This proves (\ref
{star}) as
required.
\end{pf}
\begin{pf*}{Proof of Theorem~\ref{th:inv}\textup{(ii)}}
Let us define $m_n=n(n+1)/2$. Observe that
%
%
\begin{equation}\label{zetar}\quad
\sum_{i=1}^n\rho_i = 2 S_n + D_n %
\qquad\mbox{where }
S_n:=\sum_{e \in O_n} U_e %
\quad\mbox{and}\quad D_n:=\sum_{i=1}^nU_{i,i} .
\end{equation}
Here, as in the previous proof $O_n$ denotes the set of off-diagonal edges.
For $\alpha\in(1,2)$, we have by the weak law of large numbers
$S_n/m_n \to1$
and $D_n/n \to1$ in probability. Therefore
%
%
\begin{equation}\label{as1}
\lim_{n\to\infty}\frac1{m_n}\sum_{i=1}^n \rho_i = 2
\qquad\mbox{in probability}.
\end{equation}
Theorem~\ref{th:inv}(ii) thus follows directly from Lemma~\ref{ordst} and
(\ref{as1}). The same reasoning applies in the case $\alpha=1$
replacing the
law of large numbers by the statement (\ref{eq:llnweaka1}) which now gives
(\ref{as1}) with $m_n$ replaced by $m_n w_{m_n}$.
\end{pf*}
\begin{pf*}{Proof of Theorem~\ref{th:inv}\textup{(i)}}
If $U_{e_1}\geq U_{e_2}\geq\cdots$ are the ranked values of the
i.i.d. random vector $(U_e)_{e\in O_n}$ and $S_n$ is their sum as in
(\ref{zetar}),
then by Lem\-ma~\ref{le:PoiExt}(ii), replacing $n$ with $m_n$, we have
%
%
\begin{equation}\label{pdamn}
\biggl(\frac{U_{e_1}}{S_n},\frac{U_{e_2}}{S_n},\ldots\biggr)
\mathop{\longrightarrow}^{d}_{n\to\infty}
\biggl(
\frac{x_1}{\sum_{i=1}^\infty x_i}, \frac{x_2}{\sum_{i=1}^\infty
x_i},\ldots
\biggr),
\end{equation}
where $x_1>x_2>\cdots$ denote the ranked points of the Poisson point process
on $(0,\infty)$ with intensity $\alpha x^{-\alpha-1}$.

Write $X_1,X_2,\ldots$ for the ranked values of row sums as in Lemma
\ref{ordst}, so that $\wt\rho_i = X_i/(2S_n + D_n)$, where $D_n,S_n$
are as in
(\ref{zetar}). Let
\[
Y_{2\ell-1} = \frac{X_{2\ell-1}}{2S_n + D_n} - \frac{U_{e_\ell
}}{2S_n} ,\qquad
Y_{2\ell} = \frac{X_{2\ell}}{2S_n + D_n} - \frac{U_{e_\ell
}}{2S_n} .
\]
Thanks to (\ref{pdamn}) it is sufficient to prove that
$\mathbb{P}(\max_{1\leq i\leq2k} |Y_i|> \varepsilon) \to0$, as
$n\to\infty$, for any fixed $\varepsilon>0$ and $k\in\mathbb{N}$. This
follows from the argument used in the proof of (\ref{sosh}) and
(\ref{star}).
\end{pf*}

\vspace*{-14pt}

\begin{appendix}
\section{Self-adjoint operators on PWIT}

The following classical lemma was used in Section~\ref{se:convesd}. If
$\bS$
is a self-adjoint operator on $D(\bS) \subset L^2 (V)$ with $V$ countable,
the skeleton of $\bS$ is the graph on $V$ obtained by putting an edge between
two vertices $(\bv, \bw)$ iff \mbox{$\langle\delta_\bv,\bS\delta_\bw\rangle\neq 0$}.

\begin{lem}[(Resolvent of self-adjoint operators on bipartite
graphs)]\label{le:bipartite}
Let $\bS$ be a self-adjoint operator on $D(\bS) \subset L^2 (V)$ with $V$
countable. If the skeleton is a bipartite graph then for $\bv\in V$, $h(z)
= \langle\delta_{\bv} ,(\bS- z I ) ^{-1} \delta_{\bv} \rangle$ satisfies
for all $z \in\mathbb{C}_+$, $h(-\bar z) = -\bar h (z)$.
\end{lem}
\begin{pf}
Assume first that $\bS$ is bounded: for all $\bw\in V$, $\|\bS
\delta_{\bw}\| \leq C$. For $|z| > C$, the series expansion of the resolvent
gives
\[
h(z) = - \sum_{\ell\geq0} \frac{\langle
\delta_{\bv} , \bS^\ell\delta_{\bv} \rangle}{z^{\ell+1}}.
\]
However, since the skeleton is a bipartite graph, all cycles have an even
length, and for $\ell$ odd $\langle\delta_{\bv} , \bS^\ell\delta
_{\bv}
\rangle= 0$. We deduce that for $|z| > C$, $h(-\bar z) = -\bar h (z)$. We
may then extend to $\mathbb{C}_+$ this last identity by analyticity.

If $\bS$ is not bounded, then $\bS$ is limit of a sequence of bounded
operators, and we conclude by invoking Theorem VIII.25(a) in
\cite{reedsimon}.
\end{pf}

The arguments of Section~\ref{se:convesd} were crucially based on the
following fact.

\begin{prop}\label{esssa}
The operator $\bT$ defined by (\ref{tone}) is essentially self-adjoint.
\end{prop}

To prove the proposition, we start with a deterministic lemma. Let
$V=\mathbb{N}^f$ denote the vertex set of the PWIT, and let $\cD$ be
the space
of finitely supported vectors. We write $\mathbf{u}\sim\bv$ if
$\mathbf{u}= \bv k$ or $\bv
= \mathbf{u}k$ for some $k\in\mathbb{N}$ (i.e., if $\mathbf{u},\bv
$ are neighbors) and
$\mathbf{u}\not\sim\bv$ otherwise. Let $A\dvtx\cD\to L^2(V)$ denote
the symmetric linear
operator defined by
%
%
\begin{equation}\label{asym}
\langle\delta_{\bv} , \bA\delta_{\bw} \rangle= w_{\mathbf
{u},\bv} = \overline{w}_{\bv,\mathbf{u}} ,
\end{equation}
and such that $w_{\mathbf{u},\bv}=0$ whenever $\mathbf{u}\not\sim
\bv$.
\begin{lem}[(Criterion of self-adjointness)]\label{le:criteresa}
Suppose that there exists a constant $\kappa>0$ and a sequence of connected
finite subsets $(S_n)_{n \geq1}$ in $V$, such that $ S_n \subset S_{n+1}$,
$\bigcup_n S_n = V$, and for every $n$ and $\bv\in S_n$,
\[
\sum_{\mathbf{u}\notin S_n \dvtx\mathbf{u}\sim\bv} |w_{\mathbf
{u},\bv} |^2\leq\kappa.
\]
Then the operator $\bA$ defined by (\ref{asym}) is essentially self-adjoint.
\end{lem}
\begin{pf}
It is sufficient to check that the only function $\varphi\in D(\bA^*)
\subset L^2 (V)$ such that
\[
\bA^*\varphi= \pm i \varphi
\]
is $\varphi= 0$ (see, e.g.,~\cite{reedsimon}, Theorem VIII.3). A similar
argument is used in~\cite{BLS}, Proposition 3. We deal with the case
$\bA^*
\varphi= i \varphi$, that is, for all $\mathbf{u}\in V$,
\[
i \varphi(\mathbf{u}) = \sum_{\bv\sim\mathbf{u}} w_{\mathbf
{u},\bv} \varphi(\bv).
\]
Here we use the notation $\varphi(\mathbf{u})=\langle\delta_u ,
\varphi\rangle$. Taking
conjugate, we also
have 
for all $\mathbf{u}\in V$
\[
- i \overline\varphi(\mathbf{u}) = \sum_{\bv\sim\mathbf{u}}
\overline w_{\mathbf{u},\bv} \overline\varphi(\bv) = \sum_{\bv
\sim\mathbf{u}} w_{\bv,\mathbf{u}} \overline\varphi( \bv).
\]
For any finite set $S \subset V$, we deduce
\begin{eqnarray*}
i \sum_{ \bv\in S} |\varphi( \bv) | ^2 &=& \sum_{\bv\in S }
\overline{\varphi}( \bv) ( \bA^* \varphi) (\bv) = \sum_{\bv\in
S } \overline{\varphi}( \bv) \sum_{\mathbf{u}\sim\bv} w_{\bv
,\mathbf{u}} \varphi(\mathbf{u}) \\
&=& \sum_{ \mathbf{u}\in S} \varphi(\mathbf{u}) \sum_{\bv
\sim\mathbf{u}} w_{\bv,\mathbf{u}}
\overline{\varphi} ( \bv) + \sum_{ \bv\in S } \overline{\varphi
}( \bv)
\sum_{\mathbf{u}\sim\bv\dvtx
\mathbf{u}\notin S } w_{\bv,\mathbf{u}} \varphi(\mathbf{u})\\
&&{} -
\sum_{ \mathbf{u}\in S } \varphi( \mathbf{u}) \sum_{\bv\sim
\mathbf{u}\dvtx\bv\notin S } w_{\bv,\mathbf{u}} \overline{\varphi}
(\bv) \\
&=& - i \sum_{\mathbf{u}\in S } |\varphi( \mathbf{u}) | ^2 +
\sum_{ \bv\in S }
\overline{\varphi}( \bv) \sum_{\mathbf{u}\sim\bv\dvtx \mathbf
{u}\notin S }
w_{\bv,\mathbf{u}} \varphi(\mathbf{u}) \\
&&{} - \sum_{ \mathbf{u}\in S }
\varphi( \mathbf{u}) \sum_{\bv
\sim\mathbf{u}\dvtx \bv\notin S } w_{\bv,\mathbf{u}} \overline
{\varphi} (\bv) .
\end{eqnarray*}
We obtain a Green formula,
\[
2 i \sum_{ \bv\in S} |\varphi( \bv) | ^2 = \sum_{ \bv\in S }
\overline{\varphi}( \bv) \sum_{\mathbf{u}\sim\bv\dvtx \mathbf
{u}\notin S } w_{\bv,\mathbf{u}}
\varphi(\mathbf{u}) - \sum_{ \bv\in S } \varphi( \bv) \sum
_{\mathbf{u}\sim\bv\dvtx \mathbf{u}
\notin S } \overline w_{\bv,\mathbf{u}} \overline{\varphi}
(\mathbf{u}).
\]
From Cauchy--Schwarz's inequality,
\begin{eqnarray*}
\sum_{ \bv\in S} |\varphi( \bv) | ^2 & \leq& \sum_{ \bv\in S }
|\varphi( \bv) | \sum_{\mathbf{u}\sim\bv\dvtx \mathbf{u}\notin S
} | w_{\bv,\mathbf{u}} | | \varphi(\mathbf{u})|\\
& \leq& \biggl( \sum_{ \bv\in S } |\varphi( \bv) |^2 \biggr)^{1/2}
\biggl( \sum_ { \bv\in S } \biggl( \sum_{\mathbf{u}\sim\bv\dvtx
\mathbf{u}\notin S } | w_{\bv,\mathbf{u}} | | \varphi(\mathbf{u})|
\biggr) ^2 \biggr) ^{1/2}.
\end{eqnarray*}
Now take $S = S_n$. From the assumption of the lemma, using again
Cauchy--Schwarz's inequality,
\[
\biggl( \sum_{\mathbf{u}\sim\bv\dvtx \mathbf{u}\notin S_n } | w_{\bv
,\mathbf{u}} | | \varphi
(\mathbf{u})| \biggr) ^2 \leq\kappa\sum_{\mathbf{u}\sim\bv\dvtx
\mathbf{u}\notin S_n }| \varphi
(\mathbf{u})| ^2 .
\]
Since $S_n$ is connected and the graph is a tree, if $\mathbf{u}\notin
S_n $ and
$\mathbf{u}\sim\bv$ then for any $\bv' \in S_n \setminus\bv$,
then $\mathbf{u}
\not\sim\bv' $. It follows that
\begin{eqnarray*}
\sum_{ \bv\in S_n} |\varphi( \bv) | ^2 \leq\sqrt\kappa\biggl(
\sum_{ \bv\in S_n } |\varphi( \bv) |^2 \biggr)^{1/2} \biggl( \sum_{
\mathbf{u}\in S^c_n }| \varphi(\mathbf{u})| ^2 \biggr)^{1/2}.
\end{eqnarray*}
Therefore,
\[
\sum_{ \bv\in S_n} |\varphi( \bv) | ^2 \leq\kappa\sum_{ \bv
\notin S_n }|
\varphi(\bv)| ^2.
\]
Since $\lim_n S_n = V$, as $n$ grows, the right-hand side goes to $0$, while
the left-hand side goes to $\|\varphi\|_2 ^2 $. We obtain $\varphi=
0$.
\end{pf}

Next, we need a technical lemma.
\begin{lem}\label{le:taufinite}
Let $\kappa>0$, $0 < \alpha< 2$, and let $0 < x_1 < x_2 < \cdots$ be a
Poisson process of intensity $1$ on $\mathbb{R}_+$. Define $\tau
_\kappa=
\inf\{ t \in\mathbb{N} \dvtx\sum_{k = t+1}^\infty x_k^{-2/\alpha}
\leq
\kappa\}$. Then $\mathbb{E} \tau_\kappa$ is finite and goes to $0$
as $\kappa$
goes to infinity.
\end{lem}
\begin{pf} First of all, the fact that $\tau_\kappa$ is a.s. finite
follows from the a.s. summability of $\sum_{k = 1}^\infty
x_k^{-2/\alpha}$.
We deduce also that a.s. there exists $\kappa>0$ such that $ \tau
_\kappa=
0$. From monotone convergence, it remains to check that $\mathbb{E}
\tau_\kappa<
\infty$. Let\vadjust{\goodbreak} $ n \geq1$ and $S_n =\sum_{k = 1}^\infty x_k^{-2/\alpha}
\ind_{\{x_k \geq n\}}$. From the L\'{e}vy--Khinchin formula, for
$\theta>0$,
\[
\mathbb{E} \exp( \theta S_n ) = \exp\biggl( \int_n ^\infty(
e^{\theta x^{-2 /
\alpha}} - 1) \,dx \biggr).
\]
As $n$ goes to infinity, if $\theta= o ( n^{2/\alpha} ) $,
\[
\int_n ^\infty( e^{\theta x^{-2 / \alpha}} - 1) \,dx %
\sim\frac{\theta}{2 /\alpha- 1} n^{- 2 /\alpha+ 1}.
\]
Hence, taking $\theta= (2 /\alpha- 1)n^{2 /\alpha- 1}$, we deduce
from the
Chernov bound, that for any integer $n\geq n_0$,
\[
\dP( S_n > \kappa) \leq e^{-\theta\kappa} \mathbb{E}
\exp( \theta S_n
) \leq3 e^{- c n^{2 /\alpha- 1}},
\]
where $n_0 \geq1 $ and $c = (2 /\alpha- 1) \kappa$. Also recall
(from the
Chernov bound) that if $N$ is a Poisson random variable with mean $n$,
then for
all $ t >0$,
\[
\dP( N \geq t ) \leq\exp\biggl( - t \log\frac{ t }{ n e } \biggr).
\]
Now if the event $\{\tau_\kappa> t \}$ holds, then either the number
of points
of the Poisson process $(x_k)_{k\geq1}$ in $[0,n]$ is larger than $t$
or $S_n
> \kappa$. We get for any integer $n\geq n_0$,
\[
\dP( \tau> t ) \leq e^{ - t \ln({ t }/({ n e }))} + 3 e^{- c n^{2
/\alpha- 1}}.
\]
We conclude by taking $n = \max( n_0 , t / ( 2e) ) $.
\end{pf}
\begin{pf*}{Proof of Proposition~\ref{esssa}}
We apply Lemma~\ref{le:criteresa} with $\bA$ given by $\bT$, the operator
defined by (\ref{tone}). For $\kappa>0$ and $\bv\in\mathbb{N}^f$, we
define the integer
\[
\tau_\kappa(\bv)= \inf\Biggl\{ t \geq0 \dvtx\sum_{k =
t+1}^\infty|y_{\bv k}|^{-2/\alpha} \leq\kappa\Biggr\}.
\]
The variables
$(\tau_\kappa(\bv))_\bv$ are i.i.d., and by Lemma \ref
{le:taufinite}, there
exists $\kappa>0$ such that $\mathbb{E} \tau_\kappa(\bv) < 1$. We
fix such
$\kappa$. Next, we give a green color to all vertices $\bv$ such that
$\tau_\kappa(\bv)\geq1$ and a red color otherwise. We consider an
exploration procedure starting from the root which stops at red
vertices and
goes on at green vertices. More formally, define the sub-forest $\cT
^g$ of
the PWIT where we put an edge between green vertices $\bv$ and $\bv k
$ iff
$1 \leq k \leq\tau_\kappa(\bv)$.

The sets $S_n$ appearing in Lemma~\ref{le:criteresa} are defined as follows.
If the root $\eset$ is red, we set $S_1 = \{\eset\}$. If the root is green,
we consider $ T^g_{\eset}$, the maximal subtree of $\cT^g$ that
contains the
root. It is a Galton--Watson tree with offspring distribution
$\tau_\kappa(\bv)$. Thanks to our choice of $\kappa$, $ T^g_{\eset
} $ is
almost surely finite. Let $V^g_{\eset}$ denote the set of vertices of $
T^g_{\eset} $, and consider the set $L^g_{\eset}$ of the leaves of $
T^g_{\eset} $. Note that $L^g_\eset$ is the set of vertices $\bv\in
V^g_{\eset}$ such that for all $1 \leq k \leq\tau_\kappa(\bv)$,
$\bv k$ is
red. Thus, when the root is green, we set $S_1 = V^g_{\eset} \bigcup
_{\bv\in
L^g_{\eset}} \{\bv k \dvtx1 \leq k \leq\tau_\kappa(\bv) \}$. By
construction, the set $S_1$ satisfies the condition of Lemma
\ref{le:criteresa}.\vadjust{\goodbreak}

Next, define the outer boundary of the root as $\{\eset\}$ as
$\partial\{\eset\}= \{1, \ldots, \tau_\kappa(\eset)\}$, and for
$\bv\neq\eset$, $\bv= (i_1, \ldots, i_k)$, set
\[
\partial\{\bv\} =
\{(i_1,\ldots, i_{k-1}, i_{k} +1) \} \cup\{(i_1,\ldots, i_{k},1),
\ldots,
(i_1,\ldots, i_{k},\tau_\kappa(\bv))\} .
\]
For a finite connected set $S$,
its outer boundary is defined by
\[
\partial S = \biggl( \bigcup_{\bv\in S} \partial\{\bv\} \biggr)
\Big\backslash S.
\]
To define the set $S_2$, suppose that $ \partial S_1 =\{u_1, \ldots,
u_n\}$.
The above procedure defining $S_1$ for the PWIT rooted at $\eset$ can
be now
repeated for the subtrees rooted at $u_1,\ldots,u_n$ to obtain sets
$S_1(u_1),\ldots,S_1(u_n)$. We can then define $S_2 = S_1 \cup\bigcup
_{1 \leq
i \leq n} S_1( u_i)$. Iterating this procedure, we may thus almost surely
define an increasing connected sequence $(S_n)$ of vertices with the
properties required in Lemma~\ref{le:criteresa}.\vspace*{-3pt}
\end{pf*}

\section{Tightness estimates}

Let $X$ and $K$ be the matrices defined by (\ref{levymatrix}) and
(\ref{srw}),
respectively. Recall that, when $\alpha\geq1$ we set $\kappa_n = n w_n
a_n^{-1}$, where $w_n=1$ if $\alpha>1$ and $w_n = \int_0^{a_n} x\cL
(dx)$ if
$\alpha=1$.
\begin{lem}\label{astight} %
\begin{enumerate}[(ii)]
\item[(i)]
For every $\alpha\in(0,2)$, the sequence $\mu_{a_n^{-1} X}$ is a.s. tight.
\item[(ii)]
For every $\alpha\in[1,2)$, the sequence $\mu_{\kappa_n K}$ is
a.s. tight.
\end{enumerate}
\end{lem}

We first recall a classical lemma on truncated moments and a lemma on the
eigenvalues.
\begin{lem}[(Truncated moments~\cite{Feller}, Theorem VIII.9.2)]\label{le:XV2}
For every $p > \alpha$,
\[
\mathbb{E} \bigl[ | X_{1,1} | ^p \ind_{\{| X_{1,1} | \leq t\}}
\bigr] %
\sim c (p) L(t) t^{p-\alpha},
\]
where $c(p) := \alpha/ ( p - \alpha) $. In particular, $ \mathbb{E}
[ |
X_{1,1} | ^p \ind_{\{| X_{1,1} | \leq a_n\}} ] \sim c (p)
a_n^{p} /n$.
\end{lem}
\begin{lem}[(Schatten bound~\cite{zhan}, proof of Theorem 3.32)]
If $A$ is an $n\times n$ complex Hermitian matrix then for every
$0<r\leq
2$,
%
%
\begin{equation}\label{eq:zhan}
\sum_{k=1}^n |\lambda_k(A)|^r
\leq
\sum_{i=1}^n {{\Biggl(\sum_{j=1}^n |A_{{i,j}}|^2\Biggr)}}^{r/2}.
\end{equation}
\end{lem}

\begin{pf*}{Proof of Lemma~\ref{astight}}

\textit{Proof of} (i).\quad Let us fix $r>0$. By definition of $\mu
_{X}$ we have
\[
\int_0^\infty | t | ^r \mu_{a_n ^{-1} X}(dt)
=\frac{1}{n}\sum_{k=1}^n | \lambda_k(a_n^{-1}X) |^r.\vadjust{\goodbreak}
\]
By using (\ref{eq:zhan}) we get for any $0\leq r\leq2$,
\[
\int_0^\infty | t|^r \mu_{a_n ^{-1} X}(dt)
\leq
Z_n:=\frac{1}{n}\sum_{i=1}^n Y_{n,i}
\qquad\mbox{where }
Y_{n,i}:={{\Biggl(\sum_{j=1}^na_n^{-2}|X_{{i,j}}|^2\Biggr)}}^{r/2}.
\]
We need to show that $(Z_n)_{n\geq1}$ is a.s. bounded. Assume for the
moment that
%
%
\begin{equation}\label{eq:r}
\sup_{n\geq1}\mathbb{E}(Y_{n,1}^4)<\infty
\end{equation}
for some choice of $r$. Since $Y_{n,1},\ldots,Y_{n,n}$ are i.i.d. for every
$n\geq1$, we get from (\ref{eq:r}) that
\[
\mathbb{E}\bigl( (Z_n-\mathbb{E} Z_n)^4\bigr) %
=n^{-4}\mathbb{E} {{\Biggl({{\Biggl(\sum_{i=1}^n Y_{n,i}-\mathbb{E}
Y_{n,i}\Biggr)}}^4\Biggr)}} %
=O(n^{-2}).
\]
Therefore, by the monotone convergence theorem, we get
$\mathbb{E} (\sum_{n\geq1}(Z_n-\mathbb{E} Z_n)^4)<\infty$, which gives
$\sum_{n\geq1} (Z_n-\mathbb{E} Z_n)^4<\infty$ a.s. and thus
$Z_n-\mathbb{E} Z_n\to0$ a.s. Now the sequence
$(\mathbb{E} Z_n)_{n\geq1}=(\mathbb{E} Y_{n,1})_{n\geq1}$ is
bounded by
(\ref{eq:r}), and it follows that $(Z_n)_{n\geq1}$ is a.s. bounded.

It remains to show that (\ref{eq:r}) holds, say if $0<4r<\alpha$. To this
end, let us define
\[
S_{n,a,b}:= \sum_{j=1}^n a_n^{-2}|X_{1,j}|^2
\mathbh{1}_{\{a_n^{-2}|X_{1,j}|^2 \in[a,b)\}}\qquad
\mbox{for every $a<b$.}
\]
Now $Y_{n,1}^4=(S_{n,0,\infty})^{2r}=(S_{n,0,1}+S_{n,1,\infty})^{2r}$
and thus,
%
%
\begin{equation}\label{eq:Yn}
\mathbb{E}(Y_{n,1}^4) \leq%
2^{2r-1} {{\{\mathbb{E} (S_{n,0,1}^{2r} ) + \mathbb{E}
(S_{n,1,\infty}^{2r})\}}}.
\end{equation}
We have $\sup_n \mathbb{E} (S_{n,0,1}^{2r}) < \infty$. Indeed, since
$2 r < 1$,
from the Jensen inequality,
\[
\mathbb{E} (S_{n,0,1}^{2r}) \leq(\mathbb{E} S_{n,0,1} ) ^{2r}
\]
and, by Lemma~\ref{le:XV2}, $\mathbb{E} S_{n,0,1} \sim_n \alpha/
(2- \alpha)$.

To deal with the second term of the right-hand side of (\ref{eq:Yn}), we
define
\[
M_n := \max_{1 \leq j \leq n} a_n^{-1}|X_{1,j}| \ind_{\{
a_n^{-1}|X_{1,j}| >1\}}
\]
and
\[
N_n := \#\{1\leq j\leq n\mbox{ s.t. }
a_n^{-1}|X_{1,j}| >1\}
.
\]
From the H\"{o}lder inequality, if $1/p + 1/q =1$, we have
%
%
\begin{equation}\label{eq:holderS}
\mathbb{E} (S_{n,1,\infty}^{2r}) \leq\mathbb{E} ( N_n^{2r}
M_n^{4r} ) \leq( \mathbb{E} N_n^{2rp} )^{1/p}
(\mathbb{E} M_n^{4rq} )^{1/q}.
\end{equation}
Recall that $\dP( |X_{1,2}|> a_n ) = (1 + o(1) )/ n \leq2 / n$ for large
enough $n$. Using the union bound, for large enough $n$,
\[
\dP(N_n\geq k)\leq\pmatrix{n\cr k}\dP(|X_{1,2}|>a_n)^k
\leq\frac{n^k}{k!}\frac{2^k}{n^k} = \frac{2^k}{k!}.
\]
In particular for any $\eta>0$, $\sup_n \mathbb{E} N_n ^\eta<
\infty$.
Similarly, since $L$ is slowly varying, for large enough $n$ and all
$t\geq
1$,
\[
\dP( M_n \geq t ) \leq n \dP( |X_{1,2}|> t a_n ) %
= n a_n^{-\alpha} t^{-\alpha} L( a_n t ) \leq2 t^{-\alpha}.
\]
It follows that if $\gamma< \alpha$, $\sup_n \mathbb{E} M_n^{\gamma
} < \infty$.
Taking $p$ and $q$ so that $4rq < \alpha$, we thus conclude from
(\ref{eq:holderS}) that $\sup_n \mathbb{E} (S_{n,1,\infty}^{2r}) <
\infty$.\vspace*{8pt}

\textit{Proof of} (ii).\quad Recall that for any $\alpha\in[1,2)$,
$\kappa_n = w_n n a_n ^{-1}$. Then, by using (\ref{eq:zhan}) we get
for any $0\leq r\leq2$,
\[
\int_0^\infty | t|^r \mu_{\kappa_n K}(dt) \leq
Z'_n:=\frac{1}{n}\sum_{i=1}^n \biggl(\frac{n w_n}{\rho_i}
\biggr)^{r} Y_{n,i},
\]
where
\[
Y_{n,i}:={{\Biggl(\sum
_{j=1}^na_n^{-2}|X_{{i,j}}|^2\Biggr)}}^{r/2}.
\]
From (\ref{eq:llninf}) (for $1 < \alpha< 2$) and (\ref{eq:llninfa1})
(for $\alpha= 1$), there exists $c>0$ such that a.s.,
\[
\limsup_{n \to\infty} \max_{1 \leq i \leq n} \biggl(\frac{n
w_n}{\rho_i}\biggr)^{r} < c.
\]
Hence for all $n$ large enough,
\[
Z'_n \leq\frac{c}{n}\sum_{i=1}^n Y_{n,i},
\]
and we conclude by using the same argument as in the proof of (i).
\end{pf*}
\end{appendix}

%

%
\printaddresses

\end{document}